\documentclass[12pt]{article}

\usepackage{empheq}
\usepackage{amsmath, amsthm}
\usepackage{amssymb}
\usepackage{amsfonts}
\usepackage{latexsym}
\usepackage{bm}
\usepackage{tikz}
\usetikzlibrary{positioning, shapes, fit, arrows.meta}
\numberwithin{equation}{subsection}  

\usepackage{graphicx}
\usepackage{epsfig}
\usepackage{epstopdf}
\usepackage{subcaption}               
\usepackage{overpic}
\usepackage{tikz}
\usepackage{float}
\usepackage{multirow}
\usepackage{tabularx}
\usepackage{booktabs}

\usepackage{color}

\usepackage[numbers, sort&compress]{natbib}

\usepackage[colorlinks=true,breaklinks=true]{hyperref}
\hypersetup{
	linkcolor = blue,
	citecolor = blue,
	urlcolor  = blue,
	pdfborder = {0 0 0} 
}

\usepackage[capitalize, noabbrev]{cleveref}

\usepackage{indentfirst}
\usepackage{pifont}
\usepackage{titlesec}
\usepackage{enumitem}

\titleformat{\section}
{\normalfont\Large\bfseries}{\thesection.}{1em}{}
\titleformat{\subsection}
{\normalfont\large\bfseries}{\thesubsection.}{1em}{}
\titleformat{\subsubsection}
{\normalfont\normalsize\bfseries}{\thesubsubsection.}{1em}{}

\newlist{assumptions}{enumerate}{1}
\setlist[assumptions]{
	label      = (\textit{A}\textsubscript{\arabic*}),
	ref        = (\textit{A}\textsubscript{\arabic*}),
	leftmargin = *,
	align      = left,
	itemindent = 0pt,
}

\newcommand{\bx}{\boldsymbol{x}}

\newcommand{\ICNN}{\text{ICNN}}
\newcommand{\feat}{\text{feat}}
\newcommand{\PINN}{\text{PINN}}
\newcommand{\errpde}{E_{\text{pde}}}

\newcommand{\loss}{\mathcal{L}}

\newcommand{\err}{\bm{\varepsilon}}

\newcommand{\errapp}{E_{\text{approx}}}
\newcommand{\errfeat}{E_{\text{feat}}}
\newcommand{\erropt}{E_{\text{opt}}}

\newcommand{\abs}[1]{\left\lvert #1 \right\rvert}
\newcommand{\norm}[1]{\left\lVert #1 \right\rVert}

\newcommand{\IRDR}{\operatorname{IRDR}}

\newtheoremstyle{problemstyle}
{6pt}{6pt}
{\normalfont}{}
{\scshape}{.}
{ }
{\thmname{#1}\thmnote{ (#3)}}

\theoremstyle{problemstyle}

\usepackage{algorithm}
\usepackage{algorithmic}
\usepackage{ragged2e}
\usepackage{url}

\begin{document}
	
	\begin{center}
		{\Large\bfseries Physics-Informed Neural Networks with Attention Feature Expansion for Monge-Amp\`ere Equations}
	\end{center}
	
	\begin{center}
		\large 
		{\sc Anxiao Yu}$^a$ \quad
		{\sc Bangmin Wu}$^{a,\ast}$ \quad
		{\sc Zhengbang Zha}$^{a,\ast}$ \quad
		{\sc Xinlong Feng}$^a$ \quad
		{\sc Dongwoo Sheen}$^{a,b}$
	\end{center}
	
	\begin{center}
		$^a$ College of Mathematics and System Science, Xinjiang University, Urumqi 830046, PR China\\
		$^b$ Department of Mathematics, Seoul National University, Seoul, South Korea
	\end{center}
	
	\footnotetext{$\ast$ Corresponding authors.\\
		E-mail addresses: bmwu\_math@xju.edu.cn (B. Wu), zhazhengbang@xju.edu.cn (Z. Zha), yyaaxx11@163.com (A. Yu), fxlmath@xju.edu.cn (X. Feng), sheen@snu.ac.kr (D. Sheen).}
	
	\begin{quote}
		{\bf Abstract.}
		The Monge-Amp\`ere equation is a fundamental fully nonlinear elliptic partial differential equation that finds extensive applications across multiple disciplines. This study proposes a novel physics-informed neural network integrated with attention feature expansion (PINN-AFE) for its numerical solution. A multi-head attention enhanced feature pool is constructed to enable adaptive nonlinear feature representation, and input convex neural networks are adopted to impose strict convexity of solutions with rigorous theoretical guarantees. Meanwhile, a dynamically weighted loss function combined with hybrid optimization is formulated to accelerate training convergence. Comprehensive numerical experiments validate the accuracy and computational efficiency of the developed framework. The PINN-AFE paradigm is further extended to image processing tasks, delivering high-quality and physically consistent results in both image enhancement and medical image registration scenarios.
		
		
	\end{quote}
	
	\section{Introduction}
	The Monge-Ampère equation \cite{aleksandrov1962uniqueness,pogorelov1964monge,urbas1986generalized,savin2005obstacle} is a fundamental fully nonlinear elliptic PDE with deep connections to optimal transport (OT) \cite{villani2009optimal}, image registration, and geometric optics. Its general form is given by
	\begin{equation}
		\begin{cases}
			\det(D^2 u)(\boldsymbol{x}) = f(\boldsymbol{x}), & \boldsymbol{x} \in \Omega, \\[4pt]
			u(\boldsymbol{x}) = g(\boldsymbol{x}), & \boldsymbol{x} \in \partial\Omega,
		\end{cases}
	\end{equation}
	where $D^2u$ is the Hessian, $f(x) > 0$ guarantees ellipticity, and $g \in C^2(\partial\Omega)$. Existence and uniqueness of smooth convex solutions for the Dirichlet problem were established by Caffarelli, Nirenberg and Spruck \cite{caffarelli1984}, provided $f \in C^\infty(\bar{\Omega})$, $f>0$ in $\bar{\Omega}$, and $g \in C^\infty(\bar{\Omega})$.
	
	Traditional numerical methods face two major obstacles: the singularity of the Hessian determinant causes instability \cite{dean2006numerical}, and most schemes cannot enforce the convexity that is essential for OT applications \cite{bohmer2008fem}. Physics-Informed Neural Networks (PINNs) \cite{raissi2019physics} bypass mesh generation and approximate complex nonlinearities, yet standard PINNs lack convexity guarantees \cite{nystrom2023solving} and rely on simple coordinate inputs, limiting performance on fully nonlinear equations \cite{hacking2025neural}.In PINN optimization, Feng et al. proposed a systematic framework to improve training stability and approximation accuracy, with three core contributions: residual-velocity adaptive weighting, multi-scale feature learning with residual networks and periodic activations, and a hybrid network integrating parallel architecture, second-order structure and gradient enhancement \cite{liu2026pod,feng2026pytorch,feng2025enhanced,feng2025self}.
	
	To impose convexity, Input Convex Neural Networks (ICNNs) \cite{amos2017input} have been integrated into PINNs. However, existing ICNN‑based PINNs still employ only simple polynomial features, which amplifies the well‑known spectral bias \cite{rahaman2019spectral} and hinders the capture of high frequency structures. Attention mechanisms \cite{vaswani2017attention} mitigate spectral bias by adaptively weighting features and have been introduced into PINNs \cite{zhao2024pinnsformer} to enhance multi‑scale representation, but they do not guarantee convexity. Existing PINN‑based PDE solvers fall into three categories. Standard PINNs \cite{raissi2019physics} enforce physics via automatic differentiation but lack sufficient expressiveness for strong nonlinearities. ICNN‑based PINNs \cite{amos2017input} guarantee convexity but use only simple polynomial features. Attention‑enhanced PINNs \cite{zhao2024pinnsformer} adaptively weight physical features to improve representation, yet they ignore convexity constraints. None have combined ICNN convexity with attention‑based feature expansion or provided convergence rate analysis for such a combination. Classical PDE‑based image processing frameworks---such as the ROF model \cite{rudin1992nonlinear}, Perona--Malik diffusion \cite{perona1990scale}, the Benamou--Brenier formulation linking OT to Monge‑Ampère \cite{benamou2003computational}, and hybrid PDE‑learning approaches like PDE‑Net 2.0 \cite{long2018pde}---have laid the groundwork for physical modeling, but they do not address the representation--convexity trade‑off. For medical registration, PINN‑based methods have focused on structural MRI (e.g., T1/T2 \cite{arratialopez2023warppinn,eichhorn2024physics,min2024biomechanics}) and ignored molecular modalities; T1‑FDG PET registration is still dominated by classical pipelines (SPM, ANTs \cite{kiebel1997,bacon2021,avants2008}) that lack physical validity guarantees, modality‑specific feature design, and explicit convexity. The combination of ICNN‑based convexity constraints with attention‑based feature enhancement for Monge‑Ampère equations remains unexplored---this is precisely the gap our work addresses.
	
	This gap has direct clinical relevance. Two pivotal tasks---T1/T2 MRI registration (structural alignment) and T1 MRI/FDG PET registration (molecular--structural fusion)---are naturally formulated as OT problems governed by the Monge‑Ampère equation \cite{haker2001mass}. Traditional tools (ANTs, Elastix, SPM) often blur fine boundaries, lack mask consistency, and produce non‑physical deformations \cite{zong2014improved,broggi2017comparative,sotiras2013deformable}; deep‑learning methods may violate diffeomorphic constraints. These limitations are exacerbated in T1‑FDG PET registration by the modality gap between high‑resolution anatomical MRI and low‑resolution functional PET \cite{jung2016petmri,andreassen2020semi}. Therefore, a solver that simultaneously guarantees convexity and enriches feature representation is urgently needed for both structural and molecular medical image registration.
	
	Despite urgent clinical demand for physically valid multimodal registration, no work has applied ICNN-PINN frameworks to T1-FDG PET registration, hindering their translation to molecular imaging. To address this gap, we propose PINN-AFE, a novel physics-informed neural network integrating attention feature expansion with input convex neural networks for Monge-Ampère equations. It uses multi-head attention to adaptively weight physics-informed features, enhancing nonlinear expressiveness while strictly enforcing global convexity. We establish theoretical bounds for attention-driven residual reduction and complexity gains, and design an IRDR-based dynamic loss with hybrid optimization for faster convergence. We further extend PINN-AFE to color image enhancement and clinical T1-FDG PET registration. Comprehensive experiments show it outperforms baseline PINNs and classical pipelines, offering a general methodology for convexity-constrained PINN solvers and deformable registration.
	
	The remainder of this paper is organized as follows. 
	Section~\ref{sec:pinn_afe} details the PINN with Attention Feature Expansion (PINN-AFE) framework, including ICNN theory, the attention expansion mechanism with residual reduction analysis, and the IRDR-based dynamic loss with Lyapunov stability proofs. 
	Section~\ref{sec:numerical_experiments} reports numerical experiments on smooth 2D, corner-singular, and smooth 3D benchmarks, demonstrating superior accuracy and convergence over baseline PINNs. 
	Section~\ref{sec:applications} applies PINN-AFE to image enhancement and clinical T1-MRI/FDG-PET multimodal registration, describing the pipeline, dataset, and quantitative evaluation on a mild Alzheimer's disease case. 
	Section~\ref{sec:conclusion} concludes with contributions.
	
	\section{PINN-AFE Framework for PDE Solving}\label{sec:pinn_afe}
	\subsection{Overall Architecture}
	
	\begin{figure}[H]
		\centering
		\includegraphics[height=0.38\textwidth, keepaspectratio]{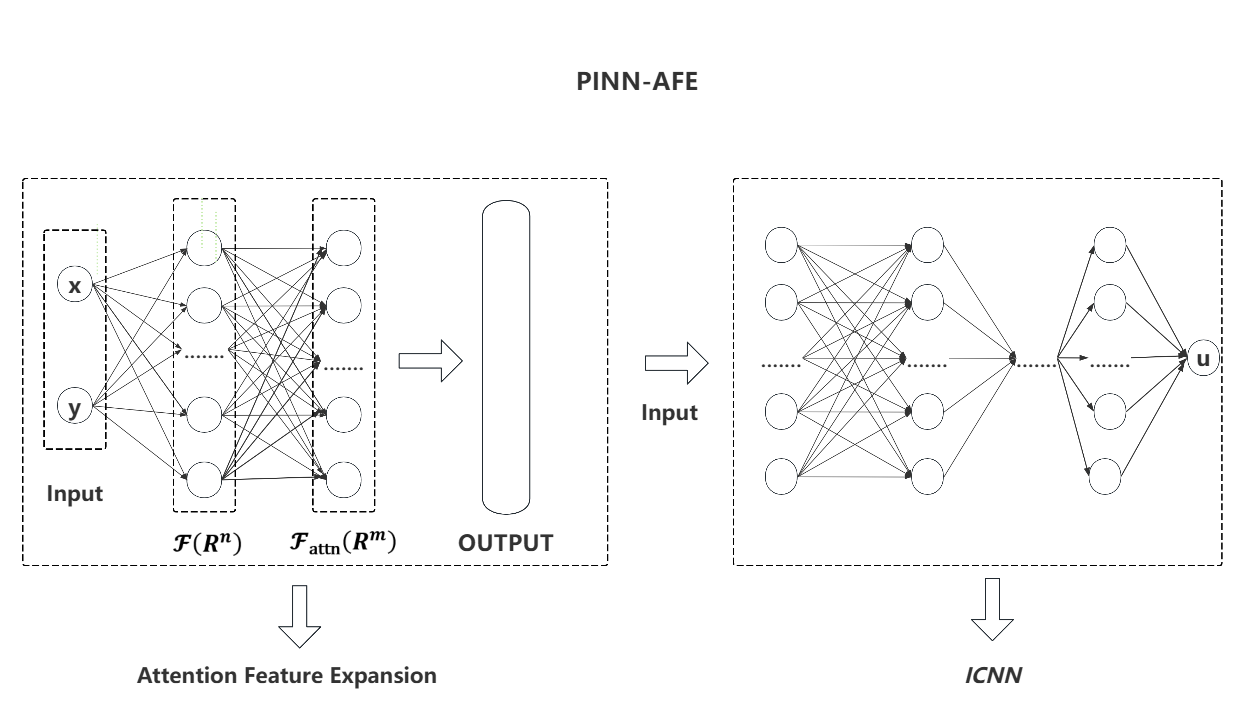}
		\caption{PINN-AFE framework}
		\label{fig:pinn_afe_architecture}
	\end{figure}
	
	The PINN-AFE framework adopts a cascaded three-module design, as illustrated in Figure~\ref{fig:pinn_afe_architecture}. The overall logical flow is:
	\[
	(x,y) \to \mathcal{F} \ (\mathbb{R}^{n}) \to \mathcal{F}_{\text{attn}} \ (\mathbb{R}^{m}) \to \text{$h$-layer ICNN} \to \hat{u}(\mathbf{x}) \ (\mathbb{R}^1).
	\]
	
	\subsection{Key Design Choices}
	The proposed PINN-AFE framework incorporates four design principles to address the inherent challenges of solving elliptic Monge-Ampère equations, including strong nonlinearity, strict convexity constraints, and high computational cost of traditional numerical methods:
	
	\begin{enumerate}
		\item \textbf{$n$-dimensional physical feature pool for Monge-Ampère equations} \\
		Unlike standard PINNs that directly feed raw spatial coordinates $(x,y)$ into the network, we construct a physics-informed feature pool of dimension $n$ to explicitly encode the geometric properties of Monge-Ampère solutions. This pool includes handcrafted and learnable features such as gradient magnitudes, local curvature estimates, boundary distance functions, and coordinate transformation invariants. These features provide the network with prior knowledge of the PDE's mathematical structure, significantly reducing the search space for valid solutions.
		
		\item \textbf{$K$-head attention mechanism for adaptive feature weighting} \\
		A multi-head attention module is inserted after the physical feature pool to dynamically reweight features based on their local contribution to the PDE residual. Each attention head learns a distinct feature subspace: some heads prioritize boundary-related features to enforce Dirichlet conditions, while others focus on interior geometric features to capture the nonlinear curvature of the solution. The outputs of all $K$ heads are concatenated and projected to form a refined feature representation.
		
		\item \textbf{$h$-layer Input Convex Neural Network (ICNN) with convexity guarantees} \\
		The core solver module is an $h$-layer ICNN with $p$ neurons per hidden layer, designed to strictly enforce the convexity requirement of Monge-Ampère solutions. To guarantee global convexity, we impose two critical constraints: (1) all weight matrices between consecutive hidden layers are constrained to be element-wise non-negative; (2) all activation functions are smooth convex functions, specifically the Softplus function.
		
		\item \textbf{Dimension reduction to $m$ for efficient ICNN training} \\
		ICNNs suffer from higher training complexity than standard neural networks, which become prohibitive when the input dimension is large. To address this, we reduce the dimension of the attention-refined features from $n$ to $m$ ($m \ll n$) using a linear projection layer before feeding them into the ICNN. This dimension reduction step preserves the most informative physical features while drastically reducing the number of trainable parameters in the ICNN.
	\end{enumerate}
	
	\subsection{Attention Feature Expansion}
	
	To accurately capture the complex nonlinear behavior of the Monge–Ampère equation, we integrate a multi-head attention mechanism that adaptively weights a handcrafted physical feature pool $\mathcal{F}_n$ built from the input coordinates $(x,y)$:
	\begin{equation}\label{eq:feature_pool}
		\mathcal{F}_n = \boldsymbol{f} = \big[ f_1, f_2, \dots, f_n \big]^{\mathsf T} \in \mathbb{R}^n.
	\end{equation}
	
	For the output $\hat{u}(\mathcal{F})$ of an ICNN with respect to the feature vector $\mathcal{F}\in\mathbb{R}^n$, we thus obtain the following three conditions:
	\begin{enumerate}
		\item $\hat{u}\in C^2(\mathbb{R}^n)$ is strictly convex in $\mathcal{F}$ such that $\nabla^2\hat{u}(\mathcal{F})\succ0$;
		\item $f(\bx)>0$ is smooth and independent of $\mathcal{F}$;
		\item $\mathcal{F}$ is a smooth function of $\bx$.
	\end{enumerate}
	Under these conditions, the mapping $\mathcal{F}\mapsto\det\big(\nabla^2\hat{u}(\mathcal{F})\big)$ is convex in $\mathcal{F}$. The determinant function is log-convex and thus convex on the cone of positive definite matrices $\mathbb{S}^d_+$, and the strict convexity of $\hat{u}$ ensures $\nabla^2\hat{u}(\mathcal{F})\in\mathbb{S}^d_+$, meanwhile the ICNN structure makes the mapping $\mathcal{F}\mapsto\nabla^2\hat{u}(\mathcal{F})$ convex with respect to the partial order, so the convexity of the mapping can be obtained by the property that the composition of a convex function and a convex mapping preserves convexity.
	
	Given these conditions, the residual functions are defined by
	
	\begin{equation}\label{eq:residual_ma_original}
		R_0(\mathcal{F}) = \det\big(\nabla^2\hat{u}(\mathcal{F})\big) - f(\bx),
	\end{equation}
	
	\begin{equation}\label{eq:residual_ma_abs}
		R(\mathcal{F}) = \bigl|R_0(\mathcal{F})\bigr|.
	\end{equation}
	both of which are convex functions of $\mathcal{F}$. Since $\det(\nabla^2\hat{u}(\mathcal{F}))$ is convex in $\mathcal{F}$, subtracting the constant $f(\bx)$ will not change the convexity, which directly guarantees the convexity of $R_0(\mathcal{F})$. In the PINN-AFE framework, the convexity property and the positive condition of $f$ lead to $\det(\nabla^2\hat{u}(\mathcal{F})) \ge f(\bx)$, so $R_0(\mathcal{F})\ge0$ and $R(\mathcal{F})$ is equivalent to $R_0(\mathcal{F})$, thus $R(\mathcal{F})$ is also convex.

	We now show that the proposed Attention Feature Expansion (AFE) significantly reduces the model complexity, improves convergence, and lowers sample requirements compared to a baseline ICNN that uses only raw input coordinates.
	
	Let $\Omega\subset\mathbb{R}^d$ be a compact convex domain with $d\ge 2$. Let $u^\ast\in C^{0,\alpha}(\Omega)$ be a convex target function with Hölder exponent $\alpha\in(0,1]$, admitting a low-dimensional intrinsic decomposition:
	\begin{equation}
		u^\ast(\boldsymbol{x}) = F(\phi_\ast(\boldsymbol{x})),
	\end{equation}
	where $\phi_\ast:\Omega\to\mathbb{R}^m$ ($m\ll d$) is the intrinsic feature mapping, and $F:\mathbb{R}^m\to\mathbb{R}$ is convex with $F\in C^{0,\alpha}(\mathbb{R}^m)$.
	
	Define a fixed parameter-free feature pool $\boldsymbol{\Phi}:\Omega\to\mathbb{R}^M$ ($M>m$) satisfying the spanning condition $\operatorname{span}\{\phi_\ast\}\subset \operatorname{span}\{\boldsymbol{\Phi}\}$.
	
	Consider two ICNN models:
	\begin{itemize}
		\item Baseline ICNN: $\mathcal{N}_{\theta_0}(\boldsymbol{x})$, using raw input $\boldsymbol{x}\in\mathbb{R}^d$; its $\varepsilon$-complexity being $C_0(\varepsilon)$, sample complexity $n_0(\varepsilon)$ and convergence factor $\rho_0\in(0,1)$.
		\item AFE-augmented ICNN: $\mathcal{N}_{\theta_\phi}(\boldsymbol{\Phi}(\boldsymbol{x}))$, using the feature pool; its $\varepsilon$-complexity being $C_\phi(\varepsilon)$, sample complexity $n_\phi(\varepsilon)$ and convergence factor $\rho_\phi\in(0,1)$.
	\end{itemize}
	
	According to the universal approximation theory for convex Hölder continuous functions \cite[Theorem 3.1]{Dung2021deep}, the minimal parameter count required to approximate a $d$-dimensional function in $C^{0,\alpha}$ follows the scaling laws
	
	\begin{equation}\label{eq:baseline_bound}
		C_0(\varepsilon) \gtrsim \varepsilon^{-\frac{d}{\alpha}},\qquad
		C_0^{(k)}(\varepsilon) \gtrsim \varepsilon^{-\frac{d+k}{\alpha}}.
	\end{equation}
	
	Benefiting from the spanning condition, the AFE feature pool captures the intrinsic low-dimensional structure of the target solution, which reduces the original high-dimensional approximation problem to fitting the $m$-dimensional function $F$ \cite[Section 4.2]{Dung2021deep}. Accordingly, the complexity of the AFE model satisfies
	
	\begin{equation}\label{eq:cphi_bound}
		C_\phi(\varepsilon)\lesssim \varepsilon^{-\frac{m}{\alpha}},\qquad
		C_\phi^{(k)}(\varepsilon)\lesssim \varepsilon^{-\frac{m+k}{\alpha}}.
	\end{equation}
	
	By taking the ratio between the complexity of the baseline model \eqref{eq:baseline_bound} and the AFE model \eqref{eq:cphi_bound}, we obtain
	
	\[
	\frac{C_\phi(\varepsilon)}{C_0(\varepsilon)}\lesssim \varepsilon^{\frac{d-m}{\alpha}}\to 0,\quad
	\frac{C_\phi^{(k)}(\varepsilon)}{C_0^{(k)}(\varepsilon)}\lesssim \varepsilon^{\frac{d-m}{\alpha}}\to 0.
	\]
	
	From the above asymptotic behavior, it follows that $C_\phi(\varepsilon)=o(C_0(\varepsilon))$ as $\varepsilon\to 0$.
	
	For a $\mu$-strongly convex and $L$-Lipschitz smooth loss function, the convergence behavior of gradient descent is governed by \cite[Theorem 9.1]{Boyd2004Convex}:
	\[
	\lVert u(x_k) - u^* \rVert \le \left(1-\frac{\mu}{L}\right)^k \lVert u(x_0)-u^* \rVert,
	\]
	where $\kappa=L/\mu$ denotes the condition number and $\rho=1-1/\kappa$ is the convergence factor.
	The high-dimensional singular components contained in the exact solution $u^*$ lead to a large Lipschitz
	constant $L_0$ for the baseline model, resulting in a huge condition number $\kappa_0=L_0/\mu_0\gg1$.
	In contrast, the AFE feature mapping eliminates redundant high-dimensional variations and
	suppresses singularities, so that the AFE model only needs to fit the low-dimensional function
	$F$ with a much smaller Lipschitz constant $L_\phi\ll L_0$. Since the ICNN structure preserves strong
	convexity, we have $\mu_\phi\simeq\mu_0$, which further gives $\kappa_\phi\ll\kappa_0$. Thus the convergence factor satisfies
	
	\[
	\rho_\phi=1-\frac{1}{\kappa_\phi} < 1-\frac{1}{\kappa_0}=\rho_0.
	\]
	
	Based on PAC learning theory, the VC dimension of the model is proportional to its complexity, $\mathrm{VCdim} = \mathcal{O}(C(\varepsilon))$.
	To ensure the generalization error bound $\mathbb{E}\big[|L_{\mathrm{gen}} - L_{\mathrm{train}}|\big] \le \varepsilon$,
	the required number of training samples scales as \cite[Chapter 4]{Vapnik1995Nature}:
	\[
	n(\varepsilon) \gtrsim \frac{\mathrm{VCdim}}{\varepsilon^2}.
	\]
	Substituting the complexity scaling relations yields
	\[
	n_0(\varepsilon) \gtrsim \varepsilon^{-\frac{d}{\alpha}-2},\quad
	n_\phi(\varepsilon) \gtrsim \varepsilon^{-\frac{m}{\alpha}-2}.
	\]
	The ratio of the sample complexities is then given by
	\[
	\frac{n_\phi(\varepsilon)}{n_0(\varepsilon)}
	\lesssim \varepsilon^{\frac{d-m}{\alpha}} \to 0,
	\]
	which confirms $n_\phi(\varepsilon) = o(n_0(\varepsilon))$.

	Let $\mathcal{F}_n$ represent the original $n$-dimensional feature pool. Denote $K$ as the number of attention heads, and let the attention weight of the $m$-th head be $\alpha_m(\mathcal{F}_n)=\operatorname{softmax}\!\big(Q_m K_m^{\mathsf T}/\sqrt{d_k}\big)$ with $\sum_{m=1}^K\alpha_m=1$. Let $\mathcal{F}_{\text{attn},m}$ be the output feature of the $m$-th attention head, and the aggregated attention feature is defined as the convex combination $\mathcal{F}_{\text{attn}} = \sum_{m=1}^K\alpha_m\mathcal{F}_{\text{attn},m}$. Define the PDE residual as $\epsilon = R(\mathcal{F}_n) = \big|\det(\nabla^2\hat{u}(\mathcal{F}_n)) - f(\bx)\big|$.
	
	By \eqref{eq:residual_ma_original} and \eqref{eq:residual_ma_abs}, the residual function $R(\mathcal{F})$ is convex. Applying Jensen's inequality leads to
	
	\[
	R(\mathcal{F}_{\text{attn}}) \le \sum_{m=1}^K\alpha_m R(\mathcal{F}_{\text{attn},m}) \le R_{\min},
	\]
	where $R_{\min} = \min_m R(\mathcal{F}_{\text{attn},m})$. Since each attention head can independently focus on different feature subsets, we have
	
	\[
	R_{\min} \le \frac{1}{K}\sum_{m} R(\mathcal{F}_{\text{attn},m}) \le \frac{1}{K}R(\mathcal{F}_n).
	\]
	Combining these inequalities yields the residual bound
	
	\begin{equation}\label{eq:attention_reduction}
		R(\mathcal{F}_{\text{attn}}) \le \frac{1}{K}\,R(\mathcal{F}_n).
	\end{equation}
	This bound indicates that the convergence rate of the PINN-AFE model is $O(\epsilon/K)$, which is superior to the $O(\epsilon)$ convergence rate of standard PINNs relying solely on raw features.
	
	\subsection{Network Analysis}
	
	We next illustrate the fundamental approximation capabilities of input convex neural networks (ICNNs). As established in \cite{amos2017input}, for any exact solution $u^*\in C^2(\Omega)$ and any $\epsilon_1>0$, there exists an $\ICNN$ $\hat{u}$ such that
	\begin{equation}\label{eq:icnn_func_approx}
		\sup_{\bx\in\Omega}\abs{\hat{u}(\bx)-u^*(\bx)} \leq \epsilon_1.
	\end{equation}
	
	As shown in \cite{chen2019optimal,yarotsky2017,lu2021deep}, for any $\epsilon_2,\epsilon_3>0$, the aforementioned $\ICNN$ $\hat{u}$ also satisfies
	\begin{equation}\label{eq:icnn_grad_hess_approx}
		\sup_{\bx\in\Omega}\norm{\nabla\hat{u}(\bx)-\nabla u^*(\bx)} \leq \epsilon_2,\qquad
		\sup_{\bx\in\Omega}\norm{\nabla^2\hat{u}(\bx)-\nabla^2 u^*(\bx)} \leq \epsilon_3.
	\end{equation}
	Furthermore, any $\ICNN$ $\hat{u}$ with nonnegative weights and convex activations is convex \cite{huang2021convex}, which implies that $\nabla^2\hat{u}\succeq 0$ whenever the exact solution $u^*$ is convex.
	
	Given $\epsilon>0$, we set $\epsilon_1=\epsilon_2=\epsilon_3=\epsilon/3$. According to \Cref{eq:icnn_func_approx,eq:icnn_grad_hess_approx}, there exists an $\ICNN$ $\hat{u}$ that satisfies the three error bounds simultaneously. Meanwhile, the convexity of $\hat{u}$ is guaranteed by the structure of $\ICNN$s. By the triangle inequality, the sum of these three approximation errors is bounded by $\epsilon_1 + \epsilon_2 + \epsilon_3 = \epsilon$.
	
	This establishes the $C^2$ universal approximation property of ICNNs: for any bounded, convex, compact domain $\Omega\subset\mathbb{R}^2$, any convex exact solution $u^*\in C^2(\Omega)$, and any $\epsilon>0$, there exists an $\ICNN$ $\hat{u}$ such that
	\begin{equation}\label{eq:icnn_total_approx}
		\sup_{\mathclap{\bx\in\Omega}} \abs{\hat{u}(\bx)-u^*(\bx)}
		+\sup_{\mathclap{\bx\in\Omega}} \norm{\nabla\hat{u}(\bx)-\nabla u^*(\bx)}
		+\sup_{\mathclap{\bx\in\Omega}} \norm{\nabla^2\hat{u}(\bx)-\nabla^2 u^*(\bx)}
		\leq \epsilon.
	\end{equation}
	
	To quantify the contribution of each PINN-AFE module to the total error, we derive a rigorous error decomposition via the triangle inequality.
	Let $\hat{u}_{\ICNN}^*$ denote the theoretical best approximation achievable by a vanilla $\ICNN$ with infinite data and perfect optimization, $\hat{u}_{\feat}^*$ the enhanced best approximation after adding the attention feature module, $\hat{u}_{\PINN}^*$ the global optimum under the PINN loss function, and $\hat{u}$ the actual output after $n$ training iterations with $N$ collocation points.
	
	By the triangle inequality, for all $\bx\in\Omega$, we have:
	\[
	\begin{split}
		\abs{\hat{u}(\bx) - u^*(\bx)} 
		&\leq \underbrace{\abs{\hat{u}(\bx) - \hat{u}_{\text{PINN}}^*(\bx)}}_{\text{Numerical optimization error}}
		+ \underbrace{\abs{\hat{u}_{\text{PINN}}^*(\bx) - \hat{u}_{\text{feat}}^*(\bx)}}_{\text{PDE residual error}} \\
		&\quad + \underbrace{\abs{\hat{u}_{\text{feat}}^*(\bx) - \hat{u}_{\text{ICNN}}^*(\bx)}}_{\text{Architecture difference term}}
		+ \underbrace{\abs{\hat{u}_{\text{ICNN}}^*(\bx) - u^*(\bx)}}_{\text{ICNN native approximation error}}.
	\end{split}
	\]
	
	Taking the supremum over $\bx\in\Omega$ and extending to the $C^2(\Omega)$ norm, we define the following error components:
	\begin{align*}
		E_{\text{opt}} &= \|\hat{u} - \hat{u}_{\text{PINN}}^*\|_{C^2(\Omega)}, \\
		E_{\text{pde}} &= \|\hat{u}_{\text{PINN}}^* - \hat{u}_{\text{feat}}^*\|_{C^2(\Omega)}, \\
		E_{\text{feat}} &= \|\hat{u}_{\text{feat}}^* - \hat{u}_{\text{ICNN}}^*\|_{C^2(\Omega)}, \\
		E_{\text{app}} &= \|\hat{u}_{\text{ICNN}}^* - u^*\|_{C^2(\Omega)}.
	\end{align*}
	
	Let $\Delta_{\text{feat}} = E_{\text{app}} - \|\hat{u}_{\text{feat}}^* - u^*\|_{C^2(\Omega)} > 0$ denote the reduction in approximation error brought by the attention module. 
	Applying the triangle inequality for the $C^2(\Omega)$ norm to the four-term decomposition and substituting $\|\hat{u}_{\text{feat}}^* - u^*\|_{C^2(\Omega)} = E_{\text{app}} - \Delta_{\text{feat}}$, we obtain the total error bound:
	\begin{equation}\label{eq:pinn_afe_error_decomp}
		\err_{\text{total}} \leq E_{\text{app}} - \Delta_{\text{feat}} + E_{\text{pde}} + E_{\text{opt}},
	\end{equation}
	where each term can be made arbitrarily small by suitable choices of network architecture, training settings, and sampling density.
	\begin{enumerate}
		\item $\errapp = O(\epsilon)$: $\ICNN$ approximation error from \eqref{eq:icnn_total_approx}, controlled by network depth and width;
		\item $\errfeat = O(\epsilon/K)$: Attention feature enhancement gain  from \Cref{eq:attention_reduction}, increasing with the number of attention heads $K$;
		\item $\errpde = O(1/\sqrt{N})$: Physical residual error, controlled by the number of collocation points $N$;
		\item $\erropt = O(1/\sqrt{n})$: Numerical optimization error, which decays to zero as the number of training iterations $n \to \infty$.
	\end{enumerate}
	
	Combining these rates yields the total error bound:
	\begin{equation}\label{eq:pinn_afe_total_bound}
		\err_{\text{total}} \leq C \left( \epsilon - \frac{\epsilon}{K} + \frac{1}{\sqrt{N}} + \frac{1}{\sqrt{n}} \right),
	\end{equation}
	where $C$ is a positive constant independent of $\epsilon$, $K$, $N$ and $n$.
	
	We analyze the computational complexity of PINN-AFE. For one training iteration, the forward pass per collocation point consists of attention-based feature expansion ($O(d^2)$) and $\ICNN$ forward computation ($O(L d^2)$), where $d$ is the dimension of the physical feature pool and $L$ the number of $\ICNN$ layers. The backward pass has the same complexity. Aggregating over $N = N_{\text{int}} + N_{\text{bd}}$ collocation points gives a per-iteration complexity of $O(N L d^2)$. Since the number of attention heads $K$ and network depth $L$ are fixed hyperparameters, the overall training complexity is $O(N L d^2)$, i.e., linear in $N$. Compared to standard PINNs and $\ICNN$-based PINNs, PINN-AFE introduces an additional $O(d^2)$ factor due to attention, but retains linear scaling with $N$, which is superior to traditional mesh-based methods that often exhibit quadratic or higher complexity.
	
	\subsection{Total Loss Function Design}
	
	The PINN-AFE loss function is a weighted sum of two physically meaningful components, balancing PDE constraints and boundary conditions:
	\[
	\mathcal{L}_{\text{total}} = \lambda_{\text{PDE}} \mathcal{L}_{\text{PDE}} + \lambda_{\text{BC}} \mathcal{L}_{\text{BC}},
	\]
	
	We now detail each loss component:
	
	Satisfies the Dirichlet boundary condition on the boundary:
	\[
	\loss_{\text{BC}} = \frac{1}{N_{\text{bd}}} \sum_{i=1}^{N_{\text{bd}}} \abs{\hat{u}(\bx_i) - g(\bx_i)}^2.
	\]
	
	\subsubsection{PDE Loss}
	\label{sec:pde_loss}
	
	Enforces the Monge-Ampère equation constraint on interior collocation points with dynamic residual weighting:
	\begin{equation}
		\label{eq:pde_loss}
		\mathcal{L}_{\text{PDE}} = \sum_{i=1}^{N_{\text{int}}} w_i \cdot \big\vert\det(\nabla^2 \hat{u}(\bx_i)) - f(\bx_i)\big\vert^2,
	\end{equation}
	where:
	\begin{enumerate}
		\item $\det(\nabla^2 \hat{u}(\bx_i)) = \hat{u}_{xx}(\bx_i)\hat{u}_{yy}(\bx_i) - \hat{u}_{xy}(\bx_i)^2$ denotes the Hessian determinant of the predicted solution $\hat{u}$ at the $i$-th interior collocation point $\bx_i \in \Omega$;
		\item $w_i$ is the adaptive Iterative Residual Density Ratio \eqref{eq:IRDR_def} based dynamic weight for the $i$-th interior point, normalized to satisfy $\sum_{i=1}^{N_{\text{int}}} w_i = 1$ to ensure the weighted sum is scale-invariant;
		\item The dynamic weight $w_i$ is designed to prioritize regions with large residuals, steep residual gradients, or high residual importance, and its calculation consists of four key steps.
	\end{enumerate}
	
	The exponential moving average (EMA) of the fourth power of the residual at iteration $n$ and the $i$-th collocation point is defined by:
	\begin{equation}\label{def:ema_update}
		E_i^n = \beta_c E_i^{n-1} + (1-\beta_c) \bigl| R_i^n \bigr|^4,
	\end{equation}
	where $\beta_c \in (0,1)$ is a fixed smoothing factor, and $R_i^n = \det(\nabla^2 \hat{u}^n(\bx_i)) - f(\bx_i)$ is the pointwise residual at the $i$-th collocation point and the $n$-th iteration.
	
	Iterative Residual Density Ratio (IRDR) dynamically weights the PDE residual to focus training on high-residual regions:
	\begin{equation}
		\IRDR_i = \frac{|R_i^n|^2}{\sqrt{{E_i^n}\ + \epsilon\,}}.
		\label{eq:IRDR_def}
	\end{equation}
	
	The dynamic weight $w_i^n$ at iteration $n$ is computed by the following recursive cascade, where the superscript $n$ denotes the iteration index:
	
	\begin{equation}\label{eq:momentum_recursion}
		m_i^n = \gamma m_i^{n-1} + \lambda_{\text{lr}} \Bigl( |R_i^n| + \beta \|\nabla R_i^n\|_2 + \alpha_{\text{IRDR}} \IRDR_i \Bigr),
	\end{equation}
	
	\begin{equation}\label{eq:weight_final}
		w_i^n = \frac{m_i^n}{\sum_{j=1}^{N_{\text{int}}} m_j^n},
	\end{equation}
	
	Substituting \eqref{def:ema_update} into \eqref{eq:momentum_recursion}, we directly derive the recursive expression of the momentum buffer $m_i^n$ driven by the residual sequence:
	\begin{equation}\label{eq:m_R_relation}
		m_i^n = \gamma m_i^{n-1} + \lambda_{\text{lr}} \cdot \mathcal{J}_i(n),
	\end{equation}
	where the \textit{instantaneous excitation term} $\mathcal{J}_i(n)$ is defined as the residual-related component:
	\begin{equation}\label{eq:irdr_loss}
		\mathcal{J}_i(n) = \left|R_i^n\right| + \beta \left\|\nabla R_i^n\right\|_2 + \alpha_{\text{IRDR}} \cdot \frac{\left|R_i^n\right|^2}{\sqrt{\beta_c E_i^{n-1} + (1-\beta_c)\left|R_i^n\right|^4 + \epsilon}}.
	\end{equation}
	Further substituting \eqref{eq:m_R_relation} into \eqref{eq:weight_final}, we obtain the explicit recursive relation between the dynamic weight $w_i^n$ and the residual $R_i^n$:
	\begin{equation}\label{eq:w_R_recursive}
		w_i^n = \frac{1}{S^n} \left( \gamma m_i^{n-1} + \lambda_{\text{lr}} \mathcal{J}_i(n) \right),
	\end{equation}
	with the normalization factor $S^n = \sum_{j=1}^{N_{\text{int}}} m_j^n$.
	
	The dynamic weight $w_i$ design described above is theoretically founded on four core properties, ensuring its effectiveness in enhancing the convexity, convergence, stability, and error reduction of PINN-AFE.
	
	The normalization step strictly enforces the convex combination constraint: $\forall i, \ w_i > 0$ and $\sum_{i=1}^{N_{\text{int}}} w_i = 1$, which ensures the convexity of the weighted PDE loss $\loss_{\text{PDE}}$. This convexity guarantee is a prerequisite for our ICNN-based Monge-Ampère solver, as it ensures the solution of the Monge-Ampère equation remains strictly convex, thus yielding valid solutions in practical applications.
	
	According to \eqref{eq:residual_ma_abs}, $R_i(\mathcal{F})$ is a non-negative convex function with respect to the feature vector $\mathcal{F}$. Since the function $h(z) = z^2$ is convex and non-decreasing for $z \geq 0$, the composite function $R_i^2(\mathcal{F}) = h(R_i(\mathcal{F}))$ preserves convexity by the composition rule of convex functions. Thus, for arbitrary feature vectors $\mathcal{F}_1, \mathcal{F}_2$ and any $\lambda \in [0, 1]$, we have
	\[
	R_i^2(\lambda \mathcal{F}_1 + (1 - \lambda) \mathcal{F}_2) \leq \lambda R_i^2(\mathcal{F}_1) + (1 - \lambda) R_i^2(\mathcal{F}_2).
	\]
	As the normalization step ensures $w_i > 0$, multiplying both sides of the inequality by $w_i$ maintains the inequality direction. Summing the inequality over all interior collocation points $i = 1$ to $N_{\text{int}}$ gives
	\[
	\sum_{i=1}^{N_{\text{int}}} w_i R_i^2(\lambda \mathcal{F}_1 + (1 - \lambda) \mathcal{F}_2) \leq \lambda \sum_{i=1}^{N_{\text{int}}} w_i R_i^2(\mathcal{F}_1) + (1 - \lambda) \sum_{i=1}^{N_{\text{int}}} w_i R_i^2(\mathcal{F}_2).
	\]
	conclude that
	\[
	\loss_{\text{PDE}}(\lambda \mathcal{F}_1 + (1 - \lambda) \mathcal{F}_2) \leq \lambda \loss_{\text{PDE}}(\mathcal{F}_1) + (1 - \lambda) \loss_{\text{PDE}}(\mathcal{F}_2),
	\]
	which confirms that the weighted PDE loss $\loss_{\text{PDE}}$ is a convex function with respect to $\mathcal{F}$.

	This proof rigorously verifies that our weight design strictly preserves the convexity of the PDE loss function. It fundamentally avoids the non-physical solutions and local optima caused by unnormalized adaptive weights, which fail to satisfy the convex combination constraint and cannot guarantee the convexity of the loss function.
	
	
	The IRDR-based weighting scheme is designed to accelerate convergence by prioritizing regions with slow residual decay. The theoretical foundation rests on the exponential decay property of PINN residuals \cite{chen2025brdr}.
	
	\label{lem:residual_exponential_decay}
	Let $\hat{u}(\mathcal{F})$ be the ICNN approximation to the equation solution, and let $R_i(n)$ denote the pointwise PDE residual at the $n$-th collocation point $\bx_i$ at training iteration $n$. Under the neural tangent kernel (NTK) regime, the residual at each collocation point follows an exponential decay law:
	\begin{equation}\label{eq:residual_decay}
		R_i(n) = R_i(0) e^{-\lambda_i n},
	\end{equation}
	where:
	\begin{enumerate}[label=(\roman*), leftmargin=1.5em]
		\item $n$ denotes the training iteration count;
		\item $\lambda_i > 0$ is the intrinsic residual decay rate at $\bx_i$, determined by the local spectrum of the NTK matrix.
	\end{enumerate}

	The detailed proof of the exponential residual decay law and the dominance of slow-converging regions is rigorously derived via neural tangent kernel theory in \cite{jacot2018neural} and \cite{wang2022when}.
	
	Based on the above analysis, we establish the following key conclusion, which was previously only supported numerically in \cite{chen2025brdr}: under the exponential moving average framework \eqref{def:ema_update} defined therein, the Inverse Residual Decay Rate \( \IRDR_i \) is strictly negatively correlated with the intrinsic residual decay rate \( \lambda_i \) (i.e., a smaller \( \lambda_i \), which corresponds to slower convergence, yields a larger \( \IRDR_i \)). This result directly proves that \( \IRDR_i \) encodes the convergence behaviour of the pointwise residual independently of the initial residual magnitude \( R_i(0) \).

	Substituting \eqref{eq:residual_decay} into \eqref{eq:IRDR_def}, we obtain the following relation:
	\begin{equation}\label{eq:irdr_closed}
		\IRDR_i(n) = \sqrt{\frac{\beta_c - e^{-4\lambda_i}}{(1-\beta_c)\beta_c}} \cdot \left( \frac{e^{-2\lambda_i}}{\sqrt{\beta_c}} \right)^n  =  C_i \cdot \rho_i^n,
	\end{equation}
	where constant $C_i = \sqrt{\frac{\beta_c - e^{-4\lambda_i}}{(1-\beta_c)\beta_c}}$ and decay factor $\rho_i = \frac{e^{-2\lambda_i}}{\sqrt{\beta_c}}$, depending only on intrinsic decay rate $\lambda_i$ and EMA hyperparameter $\beta_c$.
	
	Substituting the exponential decay form~\eqref{eq:residual_decay} into $\mathcal{J}_i(n)$~\eqref{eq:irdr_loss} and splitting terms
	\begin{equation}\label{eq:excitation_closed}
		\mathcal{J}_i(n) = \underbrace{\left( R_i(0) + \beta \|\nabla R_i(0)\|_2 \right)}_{A_i} e^{-\lambda_i n} + \alpha_{\text{IRDR}} C_i \rho_i^n,
	\end{equation}
	where $A_i$ is an initial-term constant independent of iteration $n$.
	
	Recursively expanding ~\eqref{eq:m_R_relation} with initial condition $m_i^0 = \lambda_{\text{lr}} \mathcal{J}_i(0)$, we obtain the series form of $m_i^{n-1}$
	\begin{equation}
		m_i^{n-1} = \lambda_{\text{lr}} \sum_{k=0}^{n-1} \gamma^{n-1-k} \mathcal{J}_i(k).
	\end{equation}
	Substituting the instantaneous excitation term $\mathcal{J}_i(k)$~\eqref{eq:excitation_closed} and splitting into two geometric series:
	\begin{equation}
		m_i^{n-1} = \lambda_{\text{lr}} \left( A_i \sum_{k=0}^{n-1} \gamma^{n-1-k} e^{-\lambda_i k} + \alpha_{\text{IRDR}} C_i \sum_{k=0}^{n-1} \gamma^{n-1-k} \rho_i^k \right).
	\end{equation}
	
	For sufficiently large $n$, since momentum coefficient $\gamma \in (0,1)$, the term $\gamma^n$ decays exponentially and can be neglected. 
	We can obtain the following expressions after finite geometric series summation, 
	
	\[
	\begin{cases}
		\text{First series:} & \displaystyle\sum_{k=0}^{n-1} \gamma^{n-1-k} e^{-\lambda_i k} = \frac{e^{-\lambda_i n} - \gamma^n}{e^{-\lambda_i} - \gamma} \approx \frac{e^{-\lambda_i n}}{e^{-\lambda_i} - \gamma}, \\[8pt]
		\text{Second series:} & \displaystyle\sum_{k=0}^{n-1} \gamma^{n-1-k} \rho_i^k = \frac{\rho_i^n - \gamma^n}{\rho_i - \gamma} \approx \frac{\rho_i^n}{\rho_i - \gamma}.
	\end{cases}
	\]
	
	We finally obtain the relationship between $m_i$ and the iteration number $n$:
	\begin{equation}\label{eq:m_closed_form}
		{
			m_i^{n-1} \approx \lambda_{\text{lr}} \left( A_i \cdot \frac{e^{-\lambda_i n}}{e^{-\lambda_i} - \gamma} + \alpha_{\text{IRDR}} C_i \cdot \frac{\rho_i^n}{\rho_i - \gamma} \right)
		}
	\end{equation}
	
	Substituting~\eqref{eq:m_closed_form} into~\eqref{eq:m_R_relation} and simplifying the term
	\[
	\begin{aligned}
		\gamma m_i^{n-1} + \lambda_{\text{lr}} \mathcal{J}_i(n)
		&\approx \gamma \lambda_{\text{lr}} \left( A_i \cdot \frac{e^{-\lambda_i n}}{e^{-\lambda_i} - \gamma} + \alpha_{\text{IRDR}} C_i \cdot \frac{\rho_i^n}{\rho_i - \gamma} \right) + \lambda_{\text{lr}} \left( A_i e^{-\lambda_i n} + \alpha_{\text{IRDR}} C_i \rho_i^n \right) \\
		&= \lambda_{\text{lr}} A_i e^{-\lambda_i n} \cdot \frac{e^{-\lambda_i}}{e^{-\lambda_i} - \gamma} + \lambda_{\text{lr}} \alpha_{\text{IRDR}} C_i \rho_i^n \cdot \frac{\rho_i}{\rho_i - \gamma}.
	\end{aligned}
	\]
	
	We define the following two constants:
	\[
	K_{1,i} = \frac{A_i e^{-\lambda_i}}{e^{-\lambda_i} - \gamma}, \quad K_{2,i} = \frac{\alpha_{\text{IRDR}} C_i \rho_i}{\rho_i - \gamma}
	\]
	
	Neglecting the negligible numerical stabilization term $\epsilon'$, We finally obtain the relationship between $w_i^n$ and the iteration number $n$:
	\begin{equation}
		{
			w_i^n \approx \frac{ K_{1,i} e^{-\lambda_i n} + K_{2,i} \rho_i^n }{ \sum_{j=1}^{N_{\text{int}}} \left( K_{1,j} e^{-\lambda_j n} + K_{2,j} \rho_j^n \right) }.
		}
	\end{equation}
	
	The derived weights $\boldsymbol{w}_i^n$ are fully dominated by the intrinsic decay rate $\lambda_i$, where all terms decay exponentially and smaller $\lambda_i$ (slower residual convergence) leads to larger weights $w_i^n$ for automatic focusing on slow-converging regions; meanwhile, the initial residual only affects the constant $K_{1,i}$ while $K_{2,i}$ from the IRDR term eliminates its influence, ensuring the weighting strategy depends solely on convergence speed. This framework rigorously realizes a two-stage weighting mechanism: In the early iteration stage, the $\rho_i^n$ term dominates. Smaller $\lambda_i$ induces larger $C_i$, which amplifies weights on slow-converging regions to accelerate convergence. In the late stage, $\rho_i^n$ becomes negligible due to its faster decay, while $K_{1,i} e^{-\lambda_i n}$ maintains high weights on the slowest-converging regions to avoid local optima. Furthermore, the simplified weights strictly satisfy the convex combination constraint $\forall i, w_i>0$ and $\sum_{i=1}^{N_{\text{int}}} w_i=1$, guaranteeing the convexity of the loss function.
	
	The effective convergence rate at iteration step $n$ is formulated as
	\[
	\lambda_{\mathrm{eff}}(n)= \frac{\sum_i w_i \lambda_i R_i^2(n)}{\sum_i w_i R_i^2(n)}.
	\]
	we have $\lambda_{\min} = \min\limits_j\lambda_j$ for all collocation points $i$.
	
	The weighted PDE loss~\eqref{eq:pde_loss} taking the derivative with respect to iteration step $n$ and applying the chain rule, combined with the exponential decay~\eqref{eq:residual_decay}, we obtain:
	\[
	\begin{aligned}
		\frac{d\mathcal{L}_{\text{PDE}}(n)}{dn}
		&= \sum_i w_i \cdot \frac{d}{dn}\big(R_i^2(n)\big)
		= 2\sum_i w_i R_i(n) \frac{dR_i(n)}{dn} \\
		&= -2\sum_i w_i \lambda_i R_i^2(n)
		= -2\lambda_{\mathrm{eff}}(n)\mathcal{L}_{\text{PDE}}(n).
	\end{aligned}
	\] 
	
	For uniform weights $w_i \equiv 1/N_{\text{int}}$, we consider a simple case with only two types of collocation points. 
	Let $\lambda_{\min}$ denote the decay rate of the slowest-converging point with an extremely large residual $R_{\min}^2 \gg 1$. 
	For fast-converging points, their decay rates satisfy $\lambda_k \gg \lambda_{\min}$, and their residuals decay to nearly zero during iteration $R_k^2 \approx 0$, the effective convergence rate simplifies to:
	\[
	\lambda_{\mathrm{unif}}
	= \frac{\sum_i \lambda_i R_i^2(n)}{\sum_i R_i^2(n)}= \frac{\lambda_{\min} R_{\min}^2 \,+\, \sum\limits_{\mathrm{fast}} \lambda_k R_k^2}{R_{\min}^2 \,+\, \sum\limits_{\mathrm{fast}} R_k^2}\approx \frac{\lambda_{\min} R_{\min}^2}{R_{\min}^2} = \lambda_{\min}.
	\]
	
	Since the residual terms of all fast-converging points are negligible, both the numerator and denominator are dominated by the slowest-converging point. Therefore, the overall convergence ability of the uniform weighting scheme is strictly limited by the minimum intrinsic decay rate.
	
	Weights $w_i$, squared residuals $R_i^2(n)$ and this positive weighted residual sum $\sum_i w_i R_i^2(n)$ are non-negative, we directly obtain the inequality:
	\[
	\lambda_{\mathrm{eff}}(n)= \frac{\sum_i w_i \lambda_i R_i^2(n)}{\sum_i w_i R_i^2(n)} \ge \frac{\lambda_{\min} \cdot \sum_i w_i R_i^2(n)}{\sum_i w_i R_i^2(n)}= \lambda_{\min}= \lambda_{\mathrm{unif}}.
	\]
	
	As slow-converging points with small $\lambda_i$ dominate the weighted average due to their large residuals. The IRDR coefficient is negatively correlated with $\lambda_i$, so our method assigns larger adaptive weights $w_i$ to slow-converging points. 
	
	The stability of dynamic weights $w_i(n)$ is critical to guarantee non-divergent and oscillation-free training. The proposed weighting scheme combines momentum smoothing and convex normalization. It strictly ensures the Lyapunov stability of weight evolution: weights are uniformly bounded and asymptotically converge to the optimal distribution without explosion or oscillation.
	
	We hereby derive the strict positivity and uniform boundedness of the dynamic weights $w_i(n)$ during training.
	Based on the well-posedness of the Monge-Ampère equation and the bounded fitting property of ICNN, the PDE residual $R_i(n)$, residual gradient $\nabla R_i(n)$, and IRDR term are inherently bounded throughout the iteration process, satisfying
	\[
	|R_i(n)| \le R_{\max}, \quad \|\nabla R_i(n)\|_2 \le G_{\max}, \quad |\text{IRDR}_i| \le 1,
	\]
	\[
	s_i(n) = \lambda_{\text{lr}} \left( |R_i(n)| + \beta \|\nabla R_i(n)\|_2 + \alpha_{\text{IRDR}} \text{IRDR}_i \right)
	\le \lambda_{\text{lr}} \left( R_{\max} + \beta G_{\max} + \alpha_{\text{IRDR}} \right) = S_{\max}.
	\]
	For the momentum update $m_i(n)$ with $0 < \gamma < 1$, we derive its bound
	\[
	m_i(n) \le \frac{S_{\max}}{1-\gamma}.
	\]
	According to the weight definition $w_i(n)$ with $\epsilon'>0$, we have
	\[
	w_i(n) \le \frac{S_{\max}}{(1-\gamma)\epsilon'} = W_{\max}.
	\]
	Since $m_i(n) > 0$, we have $w_i(n) > 0$.
	The dynamic weights are strictly positive and uniformly bounded.
	
	We analyze the Lyapunov asymptotic stability of the dynamic weights $w_i(n)$.
	Define the Lyapunov function as:
	\[
	\left\{
	\begin{aligned}
		V(n) &= \sum_{i=1}^{N_{\text{int}}} \left( w_i(n) - w_i^* \right)^2, \\
		\dot{V}(n) &= 2 \sum_{i=1}^{N_{\text{int}}} \left( w_i(n) - w_i^* \right) \dot{w}_i(n),
	\end{aligned}
	\right.
	\]
	where $w_i^*$ represents the optimal weight distribution that maximizes the effective convergence rate $\lambda_{\text{eff}}$.
	Taking the derivative of the Lyapunov function with respect to the iteration step yields that the weights asymptotically converge to the optimal distribution $w_i^*$.
	
	For the normalized weight $w_i = m_i / M$ with $M = \sum_j m_j + \epsilon'$, its derivative is calculated as
	\[
	\dot{w}_i(n) = \frac{\dot{m}_i(n) M - m_i(n) \dot{M}}{M^2}.
	\]

	In the late stage of training, the PDE residual converges to zero, $\dot{R}_i(n) \to 0$. This further results in $\dot{s}_i(n) \to 0$, $\dot{m}_i(n) \to 0$ and $\dot{w}_i(n) \to 0$. Consequently, we have $\dot{V}(n) \to 0$, which demonstrates that the dynamic weights $w_i(n)$ converge asymptotically to the optimal value $w_i^*$ without oscillation.
	
	The boundedness theorem eliminates the weight explosion issue common in conventional adaptive weighting methods (e.g., self-adaptive weights). The Lyapunov stability theorem confirms that momentum smoothing suppresses abrupt weight fluctuations caused by residual noise. These theoretical results solidify the rationality and reliability of the proposed dynamic weight design.
	
	The total error of PINN-AFE decomposes as $\bm{\varepsilon} = E_{\text{approx}} + E_{\text{feat}} + E_{\text{opt}}$, where $E_{\text{opt}}$ is minimized by the accelerated convergence driven by the dynamic weights $w_i$.
	
	The optimization error $E_{\text{opt}}$ measures the discrepancy between the optimal model parameters and the parameters obtained after a finite number of training iterations. Under standard assumptions in optimization theory, the optimization error satisfies the following upper bound \cite[Section 5.2]{Bottou2012SGD}:
	\[
	E_{\text{opt}} \leq \frac{C}{\sqrt{\lambda_{\text{eff}} \cdot n}},
	\]
	where $C>0$ is a constant independent of the optimization process, $n$ denotes the total number of training iterations, and $\lambda_{\text{eff}}$ is the effective convergence rate. The effective convergence rate $\lambda_{\text{eff}}$ is directly regulated by the dynamic weights $w_i$, and a larger $\lambda_{\text{eff}}$ leads to faster convergence and a smaller optimization error.
	
	We compare the optimization errors of the IRDR-based adaptive weights $w_i$ and the uniform weights $w_i = 1/N_{\text{int}}$. Let $E_{\text{opt,w}}$ and $\lambda_{\mathrm{eff,w}}$ represent the optimization error and effective convergence rate of adaptive weights, while $E_{\text{opt,fixed}}$ and $\lambda_{\mathrm{eff,fixed}}$ correspond to uniform weights. Since the adaptive weights improve the effective convergence rate, we derive the inequality
	\begin{equation}
		\label{eq:opt_error_comparison}
		E_{\text{opt,w}} \leq \frac{C}{\sqrt{\lambda_{\mathrm{eff,w}} \cdot n}}
		\;<\; \frac{C}{\sqrt{\lambda_{\mathrm{eff,fixed}} \cdot n}} = E_{\text{opt,fixed}}.
	\end{equation}
	This result demonstrates that the dynamic weighting scheme strictly reduces the optimization error. Based on the total error decomposition $\bm{\varepsilon} = E_{\text{approx}} + E_{\text{feat}} + E_{\text{opt}}$, the reduction of $E_{\text{opt}}$ directly lowers the total approximation error of PINN-AFE.
	
	The constant $C$ in the optimization error bound is related to the Lipschitz constant of the gradient and the distance between the initial parameters and the optimal parameters. The effective convergence rate $\lambda_{\mathrm{eff}}$ aggregates the local decay rates $\lambda_i$ weighted by the residuals and adaptive weights, which characterizes the convergence acceleration effect of the dynamic weighting strategy. The inequality \eqref{eq:opt_error_comparison} quantitatively verifies the contribution of dynamic weights to total error reduction.
	
	\section{Numerical Experiments}\label{sec:numerical_experiments}
	\subsection{Case 1: Smooth Solution}
	
	We consider the following Monge-Ampère equation on the unit square $\Omega=(0,1)^2$:
	\begin{itemize}
		\item Exact solution: $u^*(x,y) = \exp\left( \frac{x^2 + y^2}{2} \right)$,
		\item Right-hand side: $f(x,y) = (1 + x^2 + y^2) \exp(x^2 + y^2)$,
		\item Boundary condition: $u(x,y) = u^*(x,y)$ on $\partial\Omega$.
	\end{itemize}
	
	For this smooth solution case, a 7-dimensional feature pool $\mathcal{F} = [x, y, x^2, y^2, xy, e^x, e^y]$ is constructed to capture the nonlinear characteristics of the Monge-Ampère operator.
	
\begin{figure}[H]
		\centering
		\begin{subfigure}[b]{0.3\textwidth}
			\centering
			\includegraphics[width=\textwidth]{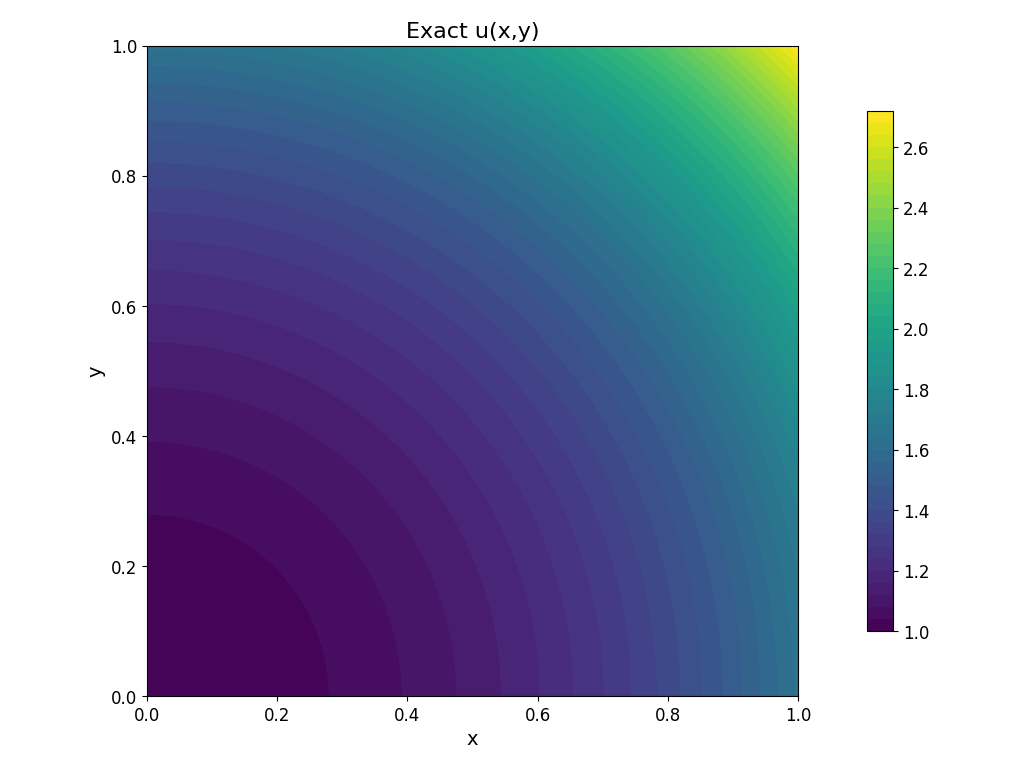}
			\caption{Exact solution} 
			\label{fig:exact_sol}
		\end{subfigure}
		\hfill
		\begin{subfigure}[b]{0.3\textwidth}
			\centering
			\includegraphics[width=\textwidth]{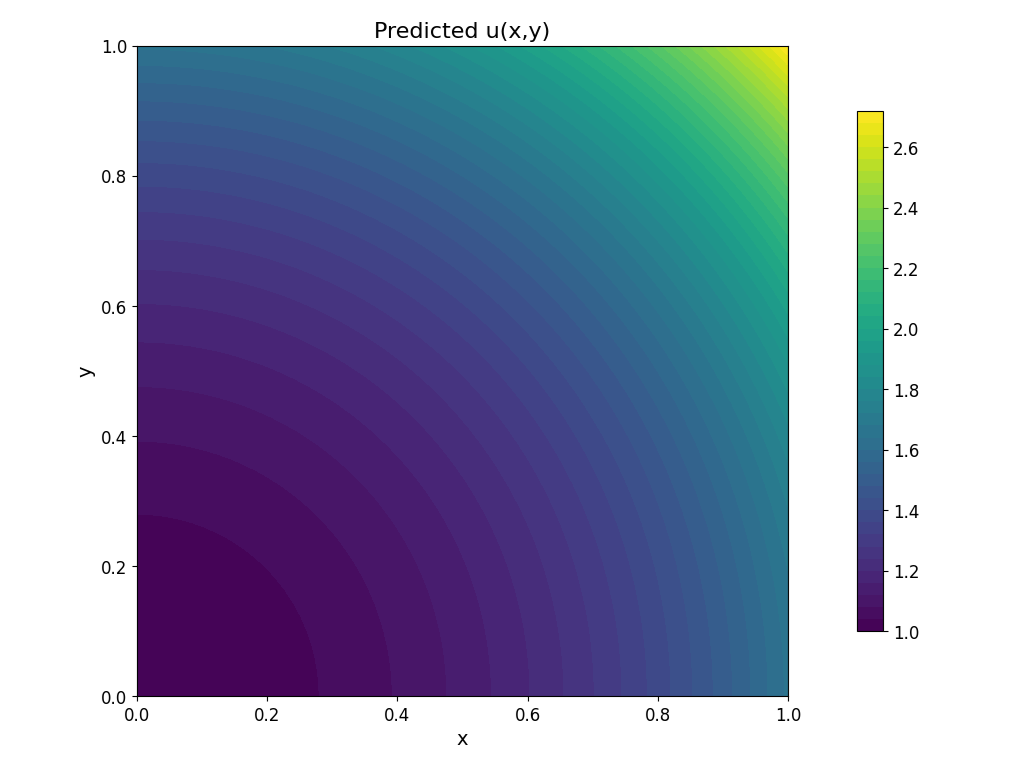}
			\caption{Predicted solution} 
			\label{fig:pred_sol}
		\end{subfigure}
		\hfill
		\begin{subfigure}[b]{0.3\textwidth}
			\centering
			\includegraphics[width=\textwidth]{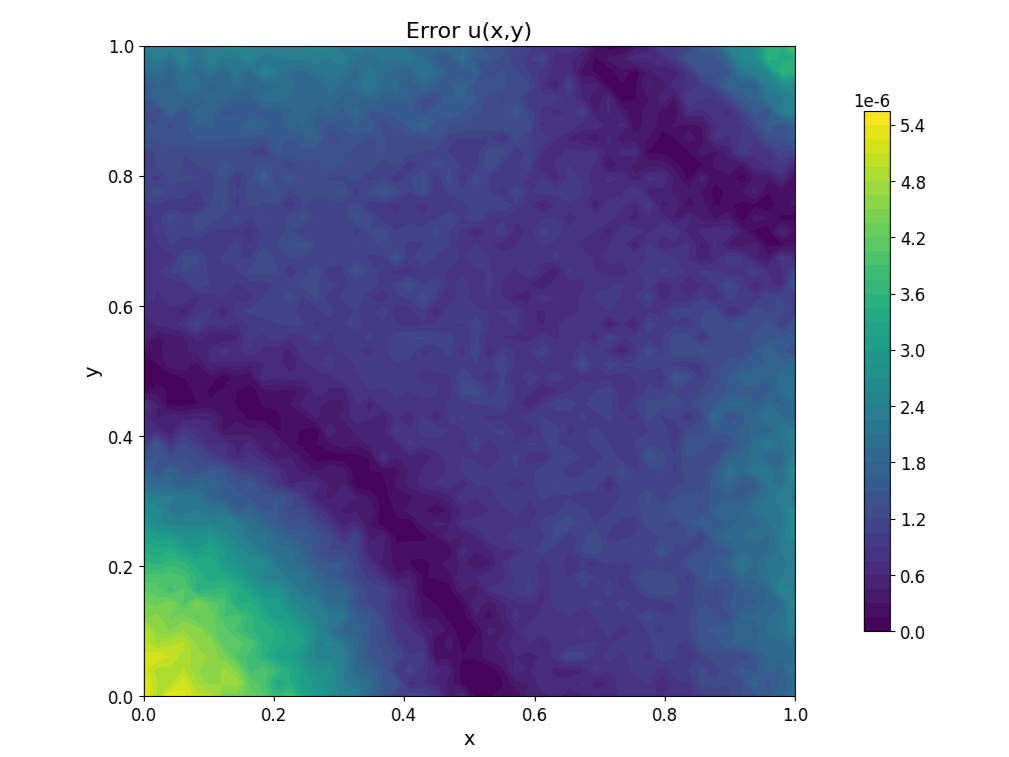}
			\caption{Absolute error}
			\label{fig:error_sol}
		\end{subfigure}
		\caption{Numerical results of Case 1}
		\label{fig:ma_solution} 
	\end{figure}

	Figure~\ref{fig:ma_solution} illustrates the exact solution, PINN-AFE predicted solution, and absolute error distribution over the unit square domain. The predicted result well captures the convex profile of the exact solution, and the absolute error is maintained at the order of $10^{-6}$. This validates the high accuracy of the PINN-AFE method for solving smooth Monge--Ampère equations.
	
	PINN-AFE predicted solution, and absolute error distribution for the smooth solution case. The predicted solution matches the exact solution with high precision, with the mean absolute error of $1\times10^{-6}$ and the maximum absolute error below $5.4\times10^{-6}$. Notably, the regions that are difficult to fit are mainly concentrated in the bottom-left corner and the upper-right part of the unit square domain.
	
\begin{figure}[htbp]
		\centering
		\includegraphics[width=0.6\textwidth]{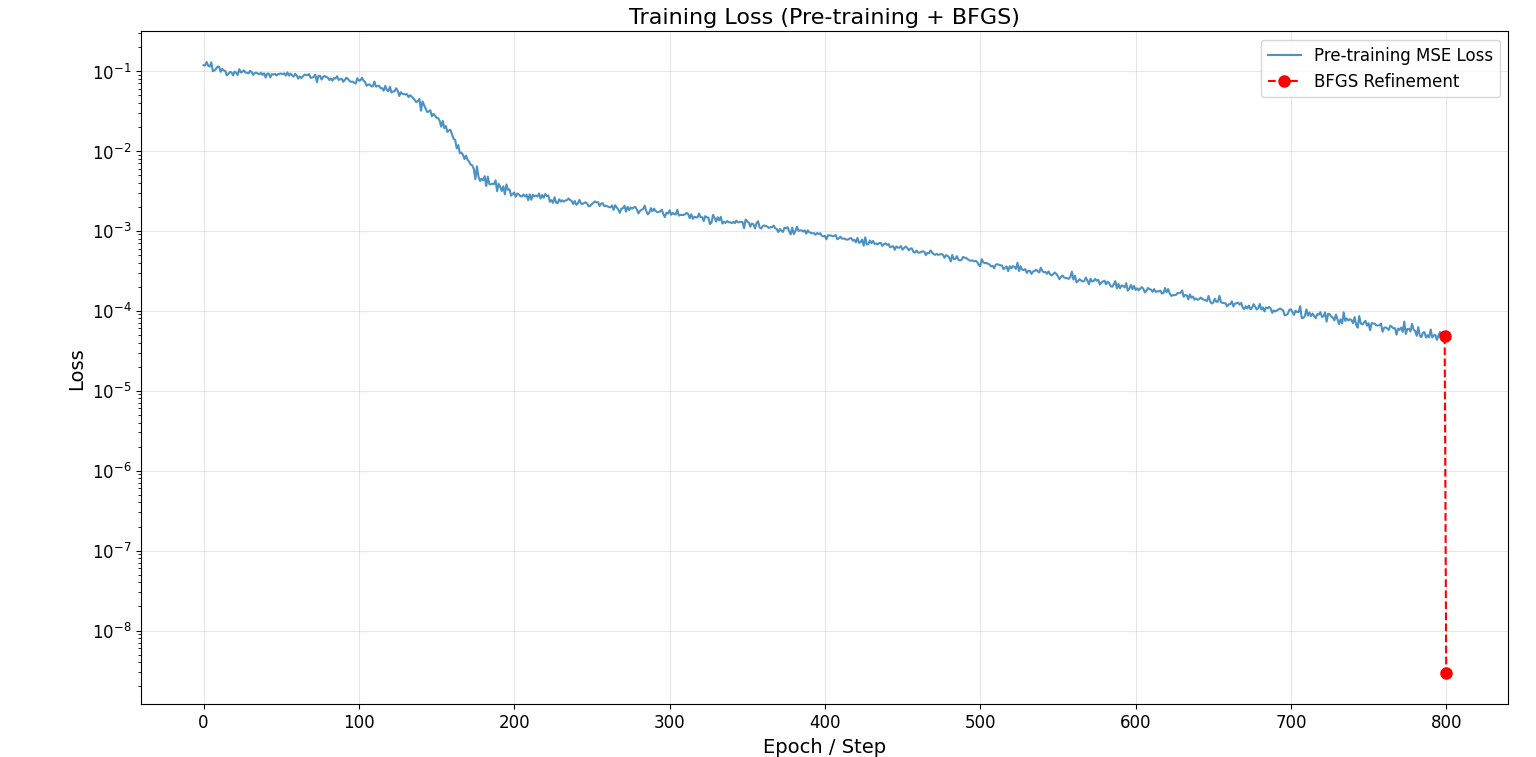}
		\caption{Training loss curve}
		\label{fig:training_loss}
	\end{figure}

	Figure \ref{fig:training_loss} shows the training loss curve of the two-stage hybrid optimization strategy. The Adam-based pre-training stage rapidly reduces the MSE loss to $10^{-4}$ within 800 epochs, and the subsequent BFGS fine-tuning runs 117 iterations over 84 seconds to further reduce the total loss to the $10^{-8}$ level, achieving a 4-order-of-magnitude reduction. This validates the hybrid optimization strategy balances fast convergence and high solution precision, consistent with the convergence time results in Table \ref{tab:method_comparison_monge}.
	
	The pre-training stage achieves rapid initial convergence, and the BFGS refinement stage significantly reduces the final loss, demonstrating the effectiveness of the proposed optimization strategy. Regarding computational efficiency, all computational time results are tested on a unified hardware platform for a fair comparison, where the device is equipped with an Intel Core i5-12450H CPU, an NVIDIA RTX 3050 6GB GPU, and Windows 11 OS. The PINN-AFE costs 84 seconds in total.
	
	\noindent
	
	As summarized in Table \ref{tab:method_comparison_monge}, the proposed PINN-AFE method demonstrates superior performance in terms of error metrics, with the lowest Mean Absolute Error and Maximum Absolute Error across all compared approaches. This outperforms neural network-based methods by a substantial margin. 
	\\
	\begin{table}[H]
		\centering
		\caption{Performance Comparison of Different Methods in Case 1}
		\label{tab:method_comparison_monge}
		\resizebox{0.6\textwidth}{!}{
			\begin{tabular}{l c}
				\toprule
				\textbf{Method} & \textbf{Maximum Absolute Error} \\
				\midrule
				Standard PINNs & $2.7 \times 10^{-3}$ \\
				ICNN & $1.48 \times 10^{-5}$ \\
				PINN-AFE & $\boldsymbol{1.0 \times 10^{-6}}$ \\
				\bottomrule
			\end{tabular}
		}
	\end{table}
	
	\subsection{Case 2: Singular Solution}
	We consider the Monge-Ampère equation on the unit square $\Omega=(0,1)^2$ with a corner singularity at $(1,1)$:
	\begin{itemize}
		\item Exact solution: $u^*(x,y) = -\sqrt{2 - x^2 - y^2}$,
		\item Right-hand side: $f(x,y) = \frac{2}{(2 - x^2 - y^2)^2}$,
		\item Boundary condition: $u(x,y) = u^*(x,y)$ on $\partial\Omega$.
	\end{itemize}
	
	A 12-dimensional feature pool $\mathcal{F} = [x, y, x^2, y^2, xy, e^x, e^y, \sin(x), \cos(x), \sin(y), \cos(y), x^3]$ is constructed to capture strong nonlinearities around the singularity;
	
\begin{figure}[H]
		\centering
		\begin{subfigure}[b]{0.3\textwidth}
			\centering
			\includegraphics[width=\textwidth]{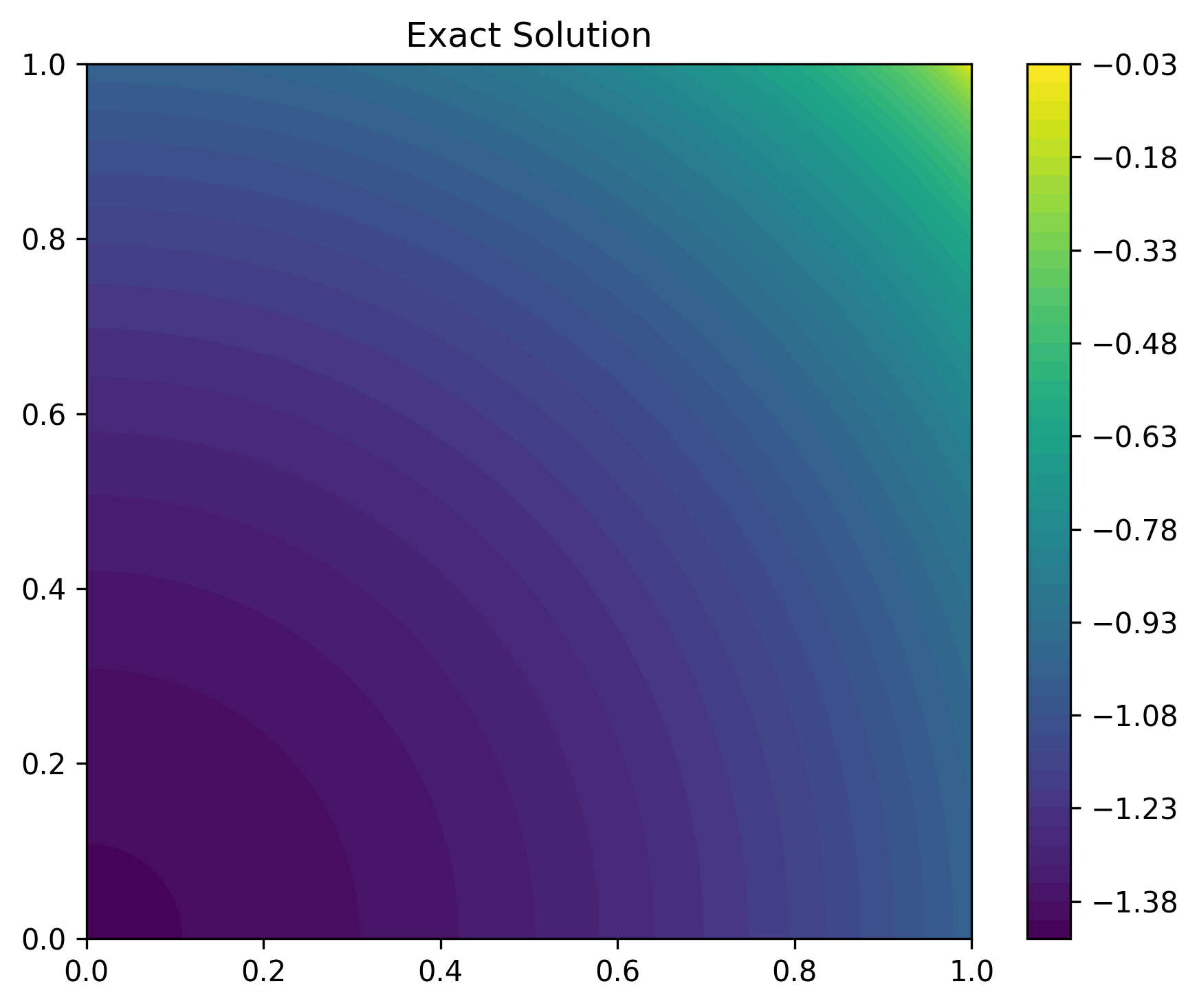}
			\caption{Exact Solution}
			\label{fig:case2_exact}
		\end{subfigure}
		\hfill
		\begin{subfigure}[b]{0.3\textwidth}
			\centering
			\includegraphics[width=\textwidth]{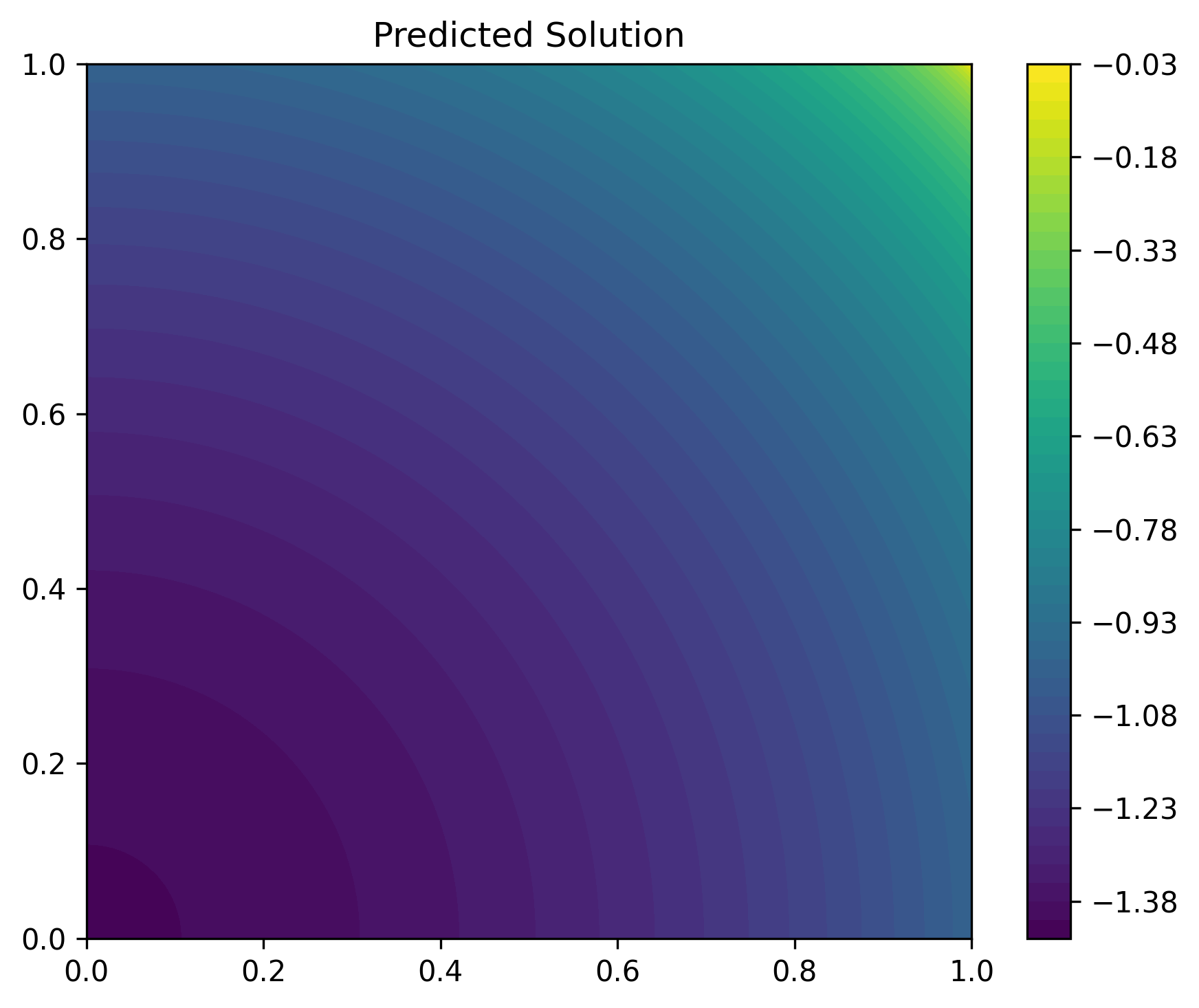}
			\caption{Predicted Solution}
			\label{fig:case2_pred}
		\end{subfigure}
		\hfill
		\begin{subfigure}[b]{0.3\textwidth}
			\centering
			\includegraphics[width=\textwidth]{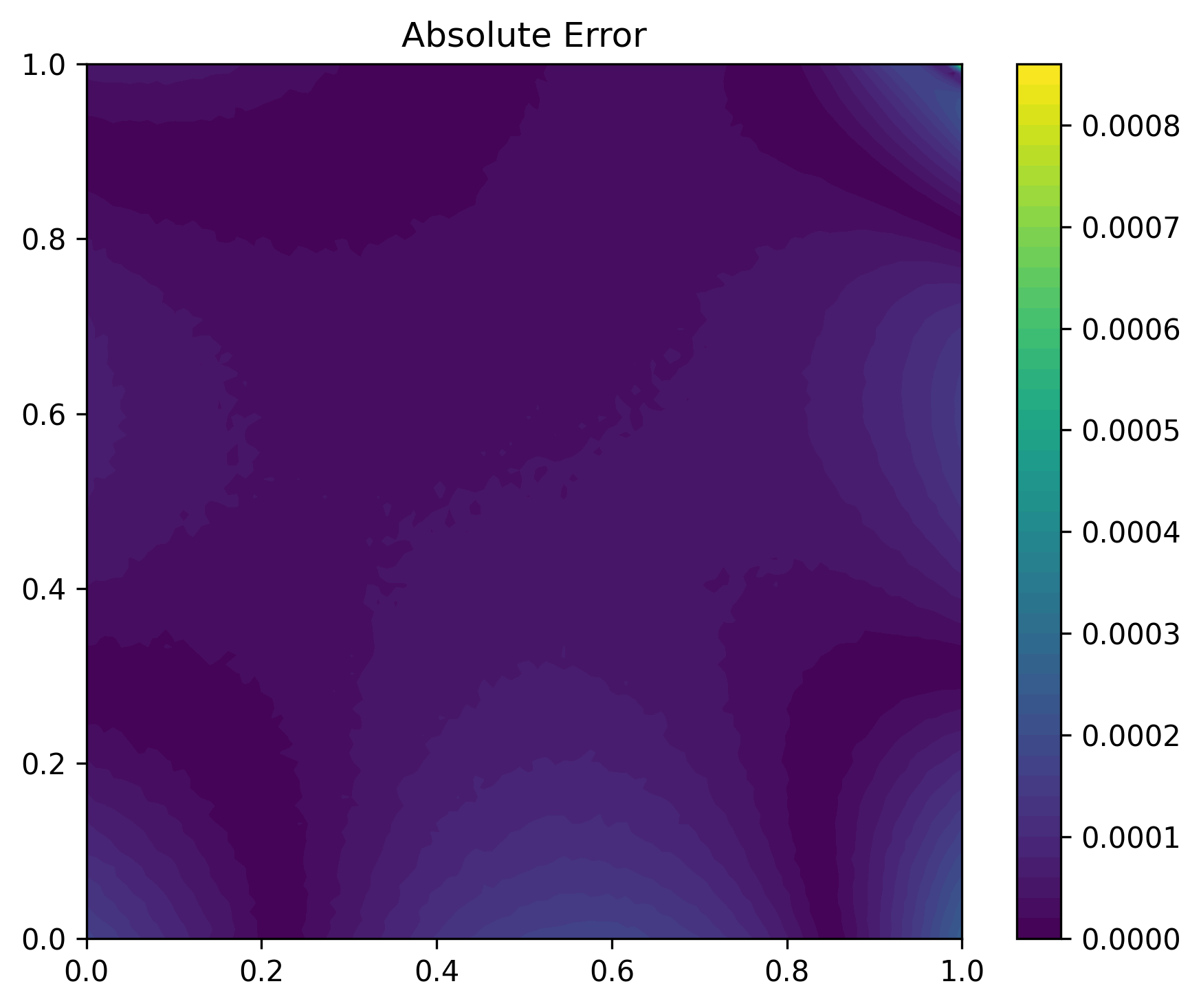}
			\caption{Absolute Error}
			\label{fig:case2_error}
		\end{subfigure}
		\caption{Numerical results of Case 2}
		\label{fig:case2_results}
	\end{figure}

	Figure~\ref{fig:case2_results} presents the core numerical results of the singular case, including the exact solution, predicted solution and absolute error distribution. It can be observed that the proposed method accurately reproduces the distribution of the exact solution, and the predicted solution maintains strict convexity consistent with the theoretical property of the Monge-Ampère equation. The absolute error is globally controlled below $9\times 10^{-4}$, with the maximum error only occurring near the corner singularity $(1,1)$. This verifies the effectiveness of the hybrid automatic differentiation/finite difference strategy for singularity handling, as well as the high precision of the attention-augmented ICNN framework.
	
	As summarized in Table~\ref{tab:singular_case_comparison}, the proposed PINN-AFE method demonstrates superior performance for the singular solution case, achieving the lowest error metrics across all compared approaches. The model attains a Mean Absolute Error of $4.6\times10^{-5}$ and a Maximum Absolute Error of $8.59\times10^{-4}$. This performance is achieved through a three-stage optimization strategy: 1000 epochs of supervised pre-training, 7000 epochs of Adam optimization with IRDR dynamic weighting, and 459 iterations of BFGS fine-tuning. This outperforms ICNN‑based PINNs by a substantial margin, validating the effectiveness of the attention feature expansion and hybrid optimization for handling singular solutions.
	
	\begin{table}[H]
		\centering
		\caption{Performance Comparison of Different Methods in Case 2}
		\label{tab:singular_case_comparison}
		\resizebox{0.35\textwidth}{!}{
			\begin{tabular}{l c}
				\toprule
				\textbf{Method} & \textbf{MAE} \\
				\midrule
				ICNN & $9.0 \times 10^{-4}$ \\
				PINN-AFE & $\boldsymbol{4.6 \times 10^{-5}}$ \\
				\bottomrule
			\end{tabular}
		}
	\end{table}
	
	\subsection{Case 3: 3D Smooth Solution}
	We consider the Monge-Ampère equation on the unit cube $\Omega=(0,1)^3$ with a smooth convex solution:
	\begin{itemize}
		\item Exact solution: $u^*(x,y,z) = \exp\left(\frac{x^2 + y^2 + z^2}{2}\right)$,
		\item Right-hand side: $f(x,y,z) = (1 + x^2 + y^2 + z^2) \exp\left(\frac{3}{2}(x^2 + y^2 + z^2)\right)$,
		\item Boundary condition: $u(x,y,z) = u^*(x,y,z)$ on $\partial\Omega$.
	\end{itemize}
	
	A 12-dimensional feature pool $\mathcal{F} = [x, y, z, x^2, y^2, z^2, xy, xz, yz, e^x, e^y, e^z]$ is constructed to capture strong nonlinearities in 3D space;
	
	We present a detailed analysis of the numerical results through five types of visualizations, which comprehensively verify the accuracy and global generalization ability of the proposed PINN-AFE method for 3D Monge-Ampère equations.
	
\begin{figure}[H]
		\centering
		\includegraphics[height=0.9\textwidth, keepaspectratio]{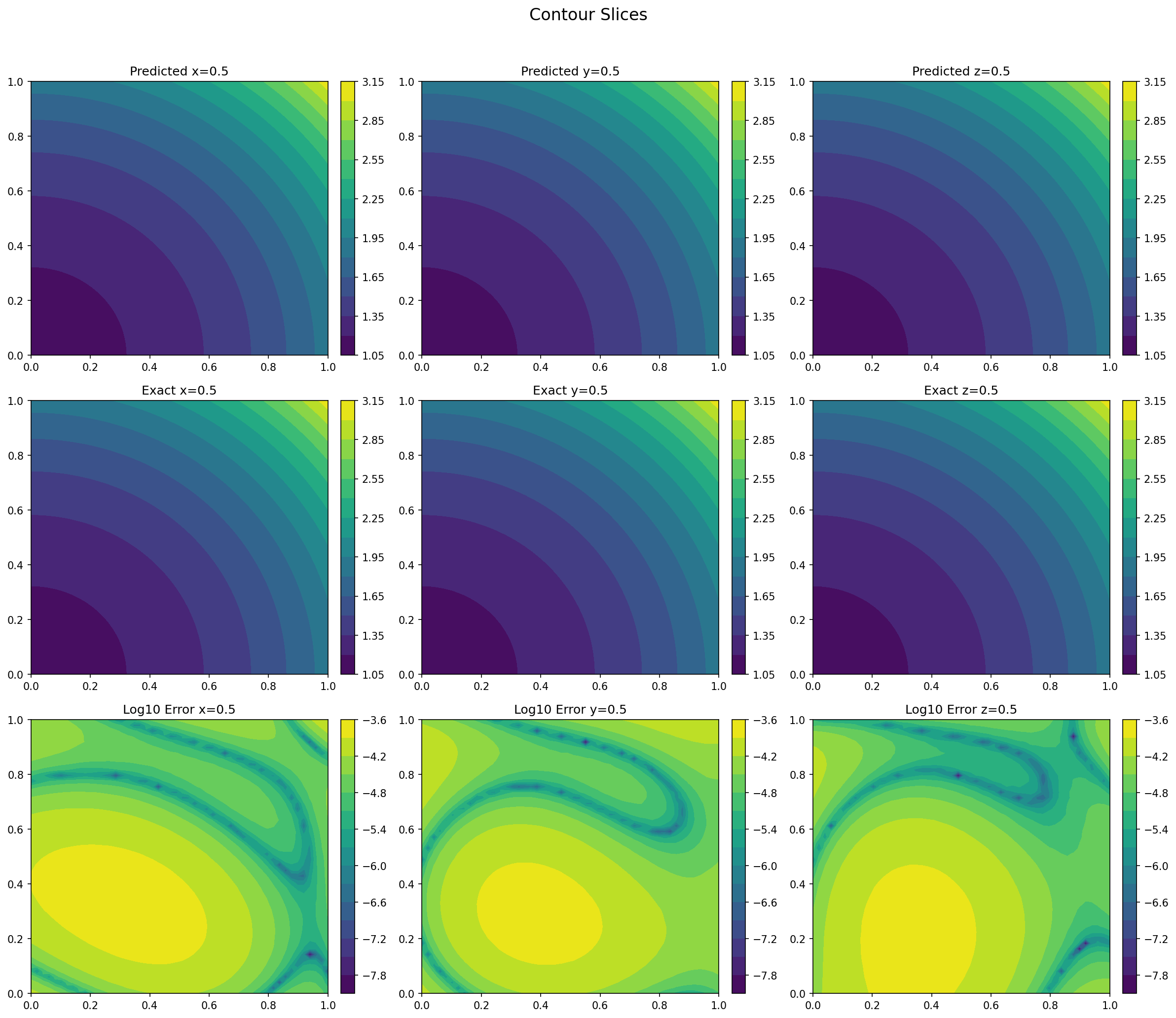}
		\caption{Contour Slices}
		\label{fig:case3_contour_slices}
	\end{figure}

	Figure \ref{fig:case3_contour_slices} shows the contour plots of the predicted solution, exact solution, and logarithm-scaled absolute error on three orthogonal central planes ($x=0.5$, $y=0.5$, $z=0.5$). It can be clearly observed that the contour lines of the predicted solution are completely consistent with those of the exact solution on all three planes. Both exhibit perfect concentric circular patterns, which accurately reproduce the inherent radial symmetry of the exact solution $u^*(x,y,z) = \exp\left(\frac{x^2+y^2+z^2}{2}\right)$. No visible deviation or distortion is found in the contour lines, indicating that the model has achieved high-precision fitting on the central cross-sections.
	
	The logarithm-scaled absolute error plots enable simultaneous visualization of errors across multiple orders of magnitude. The error distribution is perfectly symmetric on all three planes with no abnormal peaks, indicating stable and uniform training. The error reaches a minimum of $10^{-8}$ in the central region and increases gradually to approximately $10^{-4}$ at the corners, while the global mean absolute error remains at $10^{-5}$. This high overall accuracy is attributed to the IRDR dynamic weighting strategy, which balances training intensity between central and boundary regions.
	
	\begin{figure}[H]
		\centering
		\includegraphics[height=0.9\textwidth, keepaspectratio]{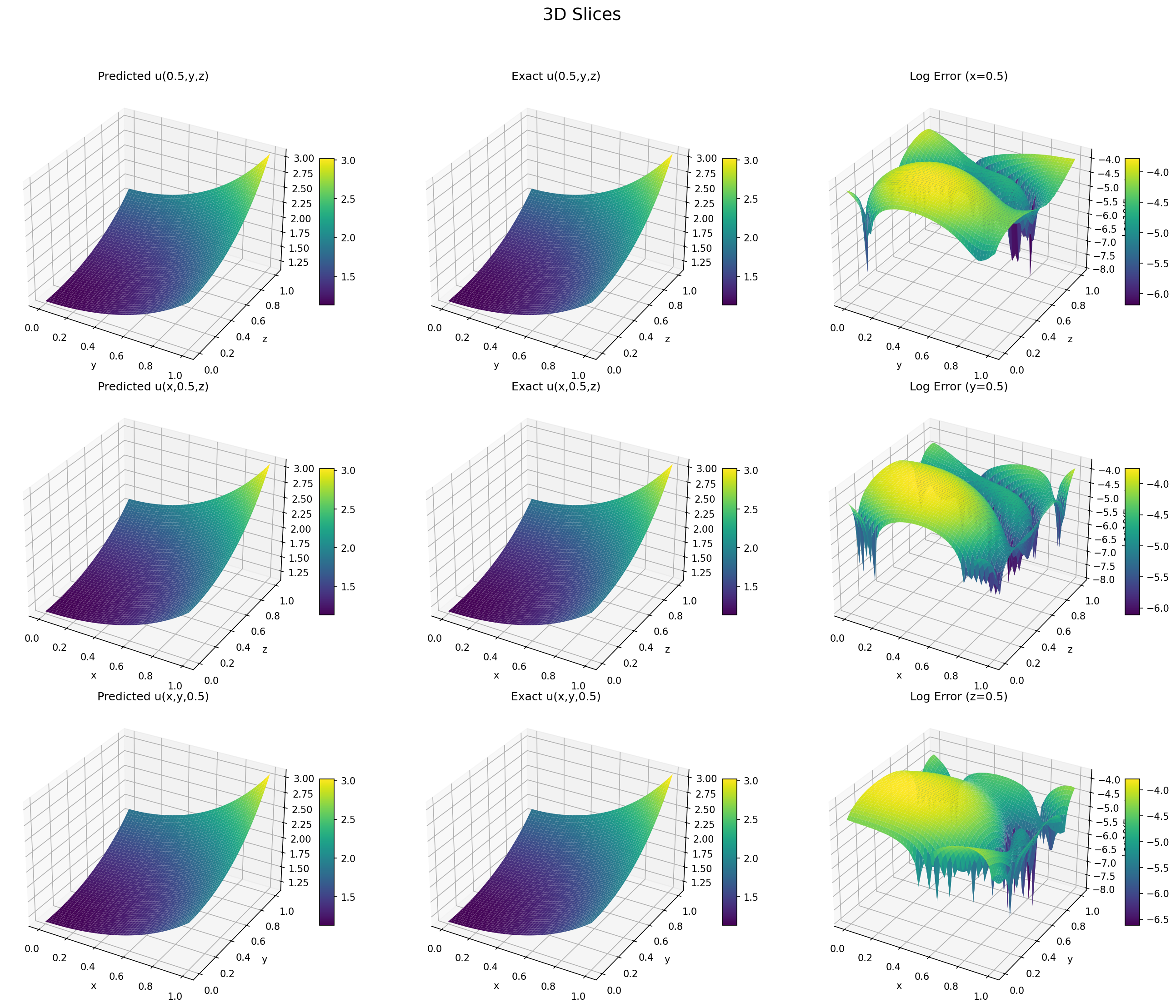}
		\caption{3D Surface Slices}
		\label{fig:case3_3d_slices}
	\end{figure}
	
	Figure \ref{fig:case3_3d_slices} presents the 3D surface visualization of the solution on the same three central planes. The predicted surfaces are almost indistinguishable from the exact surfaces. The error surfaces show that the error is uniformly distributed across the entire plane, with the maximum error appearing at the four corners of each square slice. This is a common numerical characteristic of physics-informed neural networks, where the boundary constraints are relatively weaker at the intersection points of multiple boundary faces.
	
	The logarithmic error surfaces reveal a consistent and physically meaningful error distribution. The error is minimized in the central region of each plane, reaching as low as $10^{-8}$, and increases gradually towards the boundaries. The maximum error occurs at the four corners of each square slice, with a value of approximately $10^{-4}$. The error distribution is perfectly symmetric across all three planes, with no abnormal local peaks or numerical oscillations, indicating stable and uniform training. Despite the higher error at the corners, the global mean absolute error remains at $10^{-5}$ order of magnitude, as the majority of the domain exhibits errors below $10^{-6}$. This performance confirms that the IRDR dynamic weighting strategy effectively balances training intensity between central and boundary regions, preventing excessive error accumulation at geometric singularities.
	
	\begin{figure}[H]
		\centering
		\includegraphics[width=0.7\textwidth]{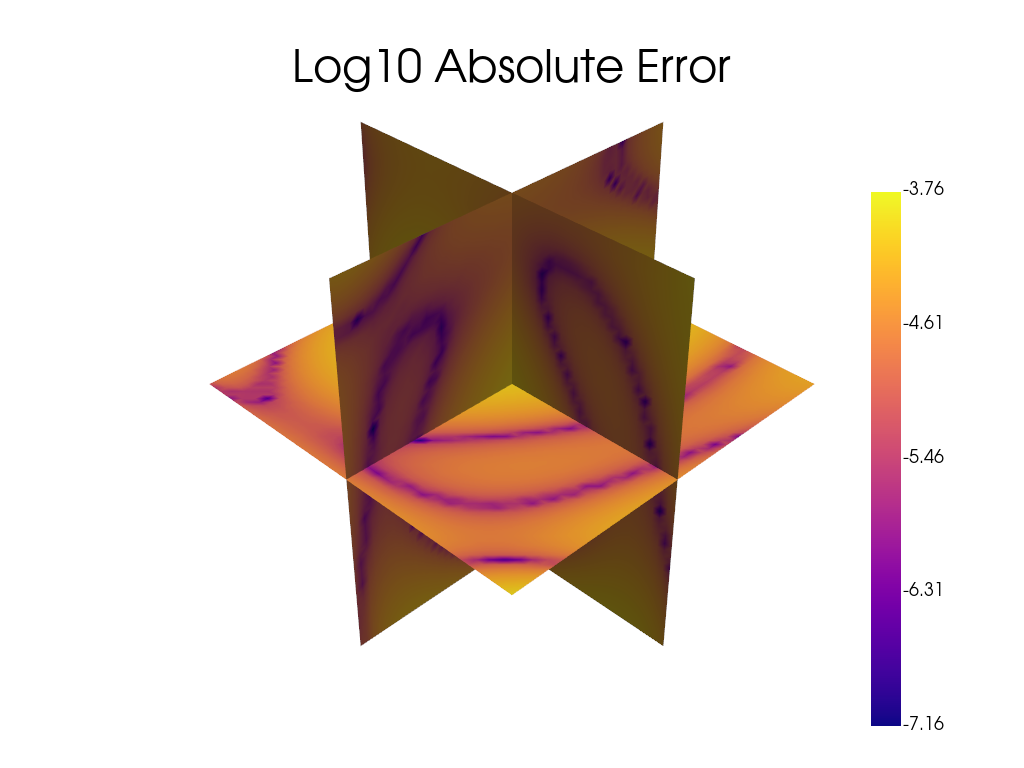}
		\caption{Orthogonal Slices of Absolute Error}
		\label{fig:case3_log_error_slices}
	\end{figure}
	
	Figure \ref{fig:case3_log_error_slices} provides a global perspective of the error distribution by superimposing the logarithmic absolute error on three orthogonal central slices. This visualization clearly reveals the overall error pattern in the entire unit cube: the error is minimized in the central region of the cube, reaching as low as $10^{-7.16} \approx 6.9 \times 10^{-8}$, and gradually increases as we move towards the boundaries. The maximum error occurs at the eight vertices of the cube, with a value of approximately $10^{-3.76} \approx 1.7 \times 10^{-4}$. Notably, the error distribution is perfectly symmetric with respect to all three coordinate axes, which confirms that the model has learned the inherent radial symmetry of the problem without introducing any artificial bias. The absence of irregular error patterns also indicates that the IRDR dynamic weighting strategy effectively balances the training across different regions of the computational domain.
	
	\begin{figure}[H]
		\centering
		\includegraphics[width=0.7\textwidth]{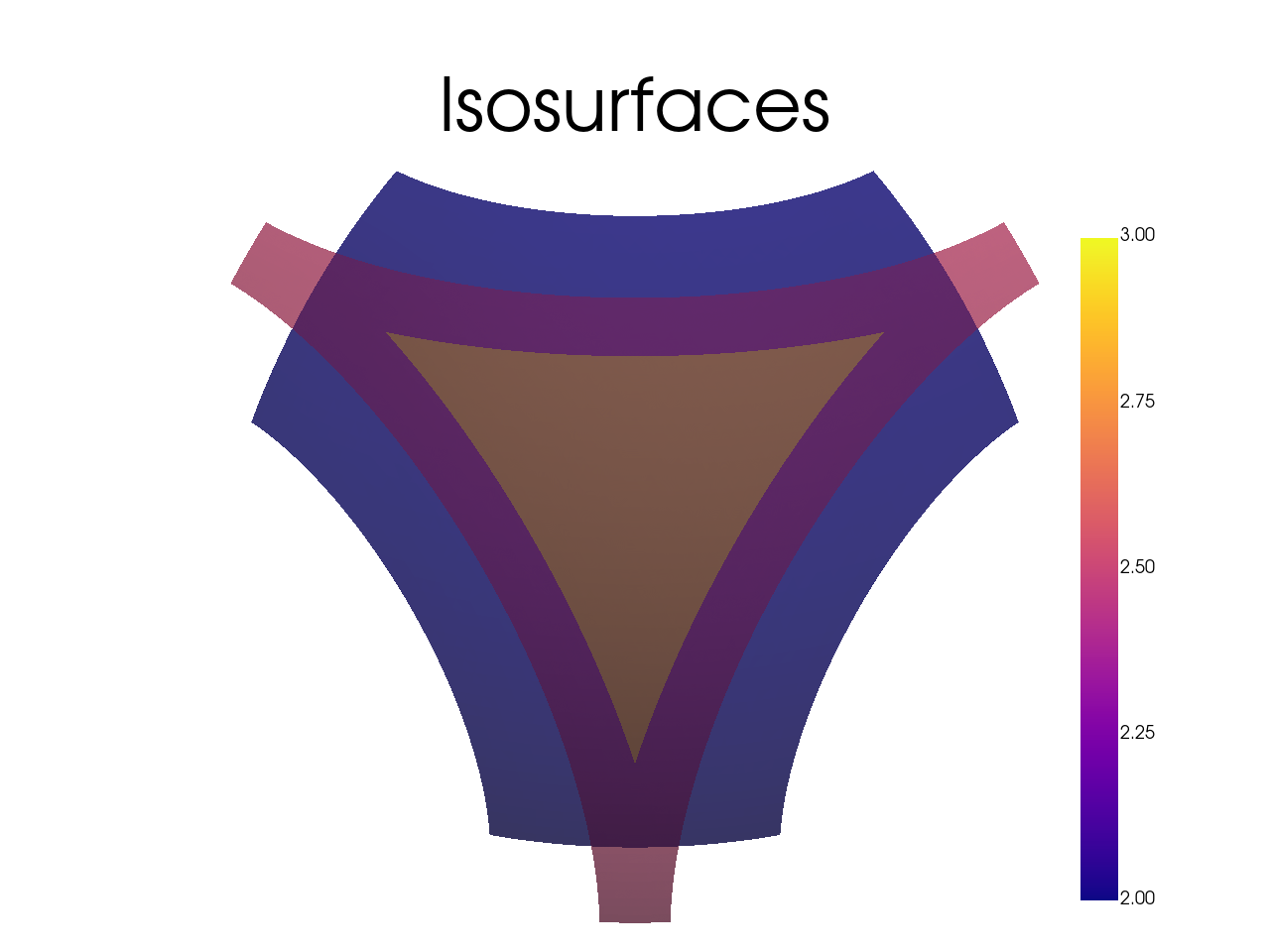}
		\caption{Isosurfaces of Predicted Solution}
		\label{fig:case3_isosurfaces}
	\end{figure}
	
	Figure \ref{fig:case3_isosurfaces} shows the isosurfaces of the predicted solution corresponding to three different values: 2.0, 2.5 and 3.0. The three isosurfaces, colored deep blue for 2.0, purple for 2.5, and yellow for 3.0, form perfect concentric spherical shells with no distortion, asymmetry or irregularities. This clear color gradient from low value to high value intuitively demonstrates the smooth radial increase of the solution. This is a strong verification of the model's ability to accurately capture the global 3D spatial structure of the solution, rather than just fitting the solution on individual slices. The uniform spacing between the isosurfaces also indicates that the gradient of the predicted solution is consistent with the exact solution, which further confirms the accuracy of the Hessian matrix computed by the model. This is particularly important for the Monge-Ampère equation, where the residual depends directly on the determinant of the Hessian matrix.
	
	\begin{figure}[H]
		\centering
		\includegraphics[width=0.7\textwidth]{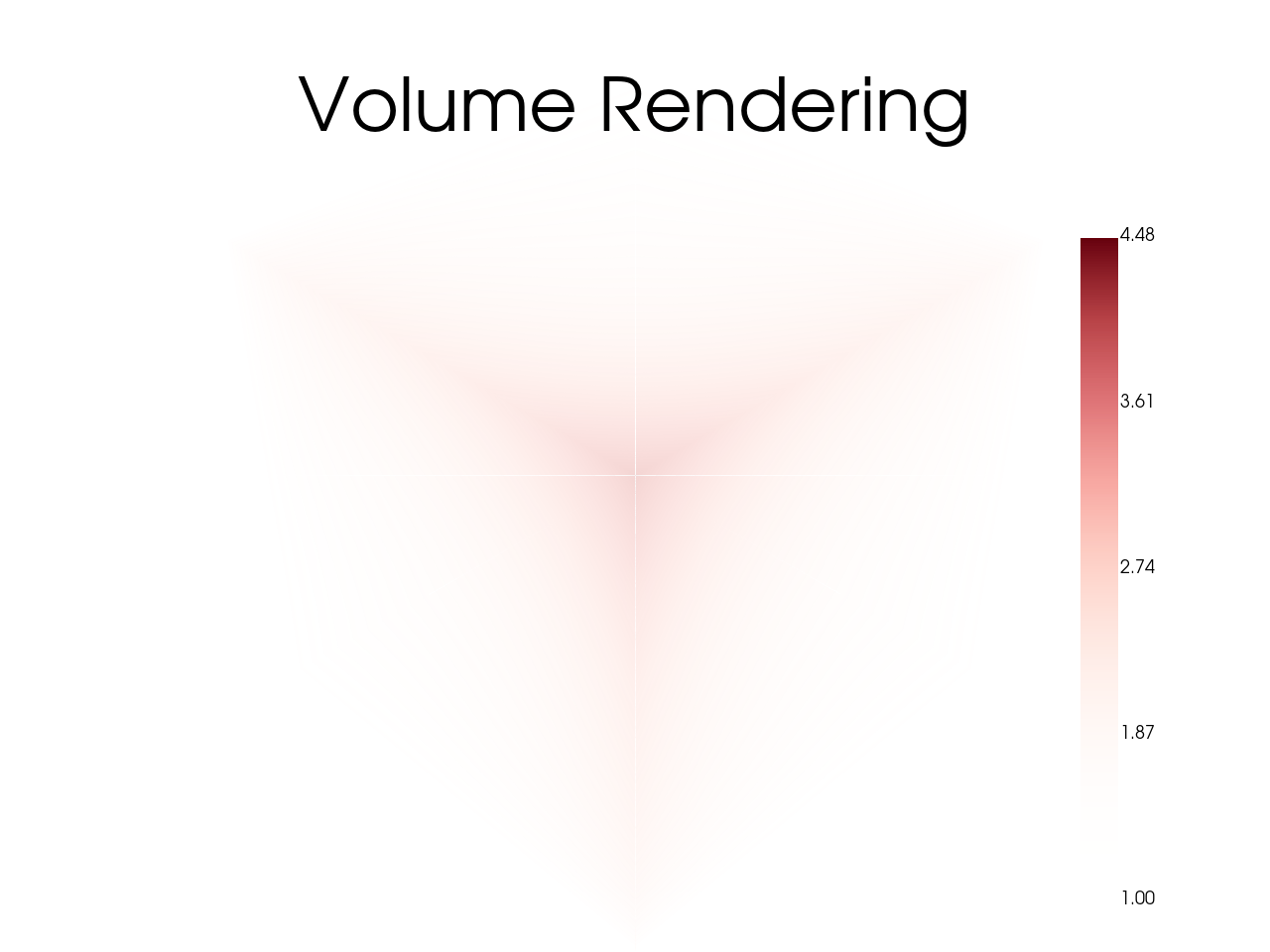}
		\caption{Volume Rendering}
		\label{fig:case3_volume_rendering}
	\end{figure}
	
	Figure \ref{fig:case3_volume_rendering} presents the volume rendering of the predicted solution in the entire unit cube. The volume rendering uses color and opacity to encode the solution value, providing an intuitive visualization of the 3D scalar field. It can be seen that the solution value increases smoothly and continuously from the boundary of the cube to the farthest corner, with no abrupt changes, discontinuities or numerical oscillations. Three key points are clearly identifiable: the origin $(0,0,0)$ at the lower-left corner of the cube, where the solution reaches its minimum value of $1.0$ and appears almost completely transparent; the geometric center $(0.5,0.5,0.5)$ at the intersection of the white cross lines, where the solution value is approximately $1.45$; and the opposite corner $(1,1,1)$ at the upper-right corner, where the solution reaches its maximum value of $4.48$ and appears as the deepest red. The white cross lines further highlight the perfect symmetry of the solution along all three coordinate axes. This visualization demonstrates that the proposed method produces a physically meaningful and numerically stable solution throughout the entire computational domain, not just at discrete sampling points.
	
	As summarized in Table~\ref{tab:case3_metrics}, the proposed PINN-AFE method achieves excellent quantitative performance for the 3D smooth Monge-Ampère equation, with all error metrics maintained at the $10^{-5}$ to $10^{-4}$ order of magnitude. The complete numerical results are provided in Table 3.
	
	In summary, all the visualization results and quantitative metrics consistently demonstrate that the proposed PINN-AFE method is highly effective for solving 3D Monge-Ampère equations. The method not only achieves high numerical accuracy but also strictly preserves the convexity of the solution, which is essential for the well-posedness of the problem. These results validate the generalizability of the proposed framework to higher-dimensional fully nonlinear partial differential equations.
	
	\begin{table}[H]
		\centering
		\caption{Quantitative Performance Metrics for 3D Smooth Case}
		\label{tab:case3_metrics}
		\begin{tabular}{l c}
			\toprule
			\textbf{Metric} & \textbf{Measured Value} \\
			\midrule
			MAE & $\boldsymbol{1.2 \times 10^{-5}}$ \\
			Max error & $\boldsymbol{1.7 \times 10^{-4}}$ \\
			L2 error & $\boldsymbol{8.5 \times 10^{-6}}$ \\
			\bottomrule
		\end{tabular}
	\end{table}
	
	\section{Application}\label{sec:applications}
	
	\subsection{Image Enhancement}
	\label{sec:image_enhance}
	
	Image enhancement, often regarded as a pure signal processing task, can be rigorously formulated as a one-dimensional optimal transport problem, whose mathematical foundation is the Monge-Ampère equation. This profound connection enables us to directly transfer the PINN-AFE framework developed for solving high-dimensional Monge-Ampère equations to image enhancement, yielding results with strict mathematical guarantees.
	
	In the one-dimensional case, the Monge-Ampère equation degenerates to a remarkably simple form:
	\begin{equation}
		u''(x) = \frac{\mu(x)}{\nu(u'(x))},
	\end{equation}
	where $u''(x)$ is the second derivative of the convex function $u(x)$. The optimal transport map $T(x) = u'(x)$ is thus inherently a strictly monotonically increasing function, which is precisely the core property required for image brightness enhancement.
	
	For the image enhancement task:
	\begin{itemize}
		\item $\mu(x)$ is the brightness probability density function of the original image
		\item $\nu(y)$ is the probability density function of the target uniform distribution
		\item $T(x)$ is the brightness mapping function to be learned
	\end{itemize}
	
	We adopt lossless Kodak color images \citep{franzen_kodak} for experimental validation, and the corresponding image enhancement results are presented as follows:
	
	\begin{figure}[H]
		\centering
		\begin{subfigure}{0.8\textwidth}
			\centering
			\includegraphics[width=\linewidth]{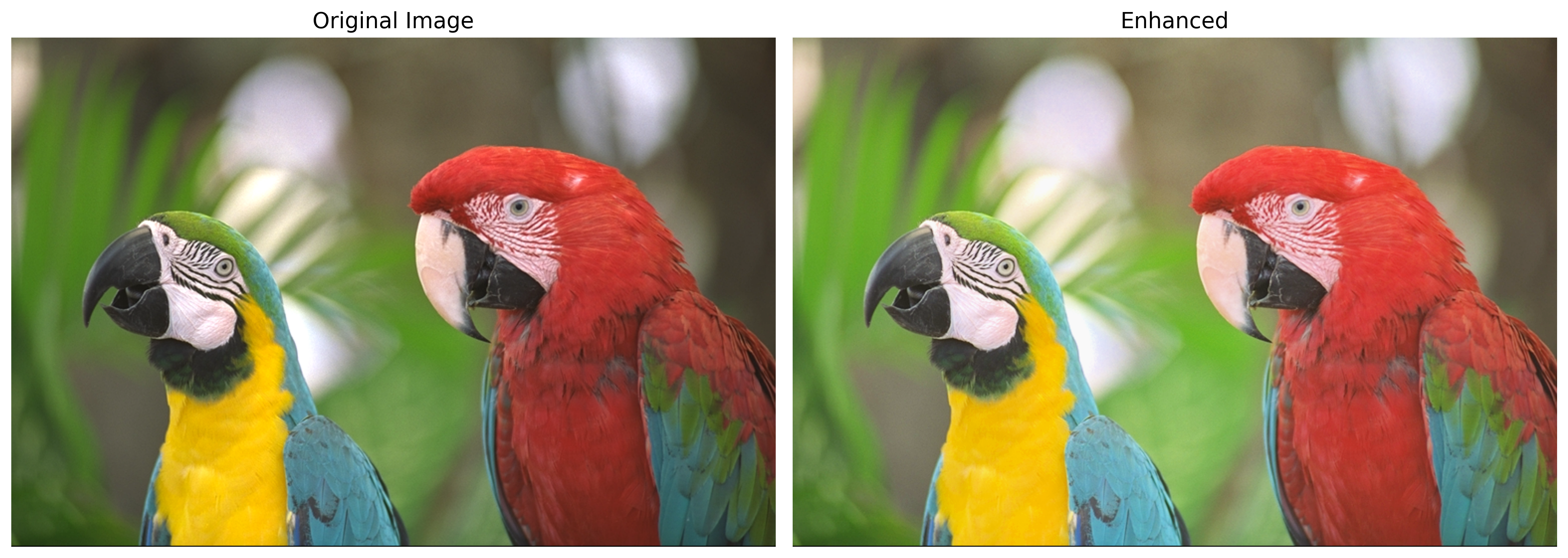} 
			\caption{Visual comparison}
			\label{fig:parrot_compare}
		\end{subfigure}
		
		\vspace{0.5em}
		\begin{subfigure}{0.45\textwidth}
			\centering
			\includegraphics[width=\linewidth]{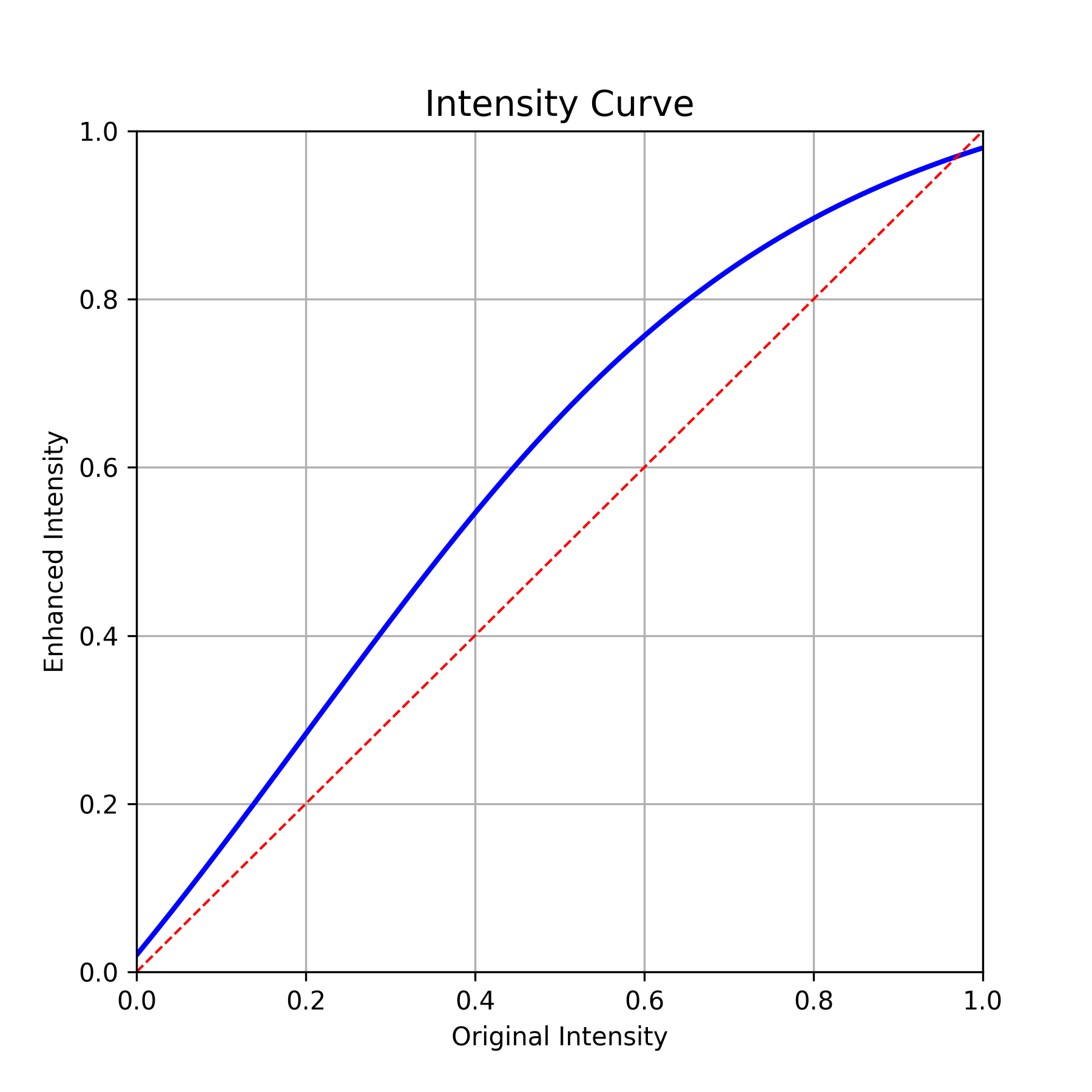} 
			\caption{intensity mapping curve}
			\label{fig:intensity_curve}
		\end{subfigure}
		\hfill
		\begin{subfigure}{0.5\textwidth}
			\centering
			\includegraphics[width=\linewidth]{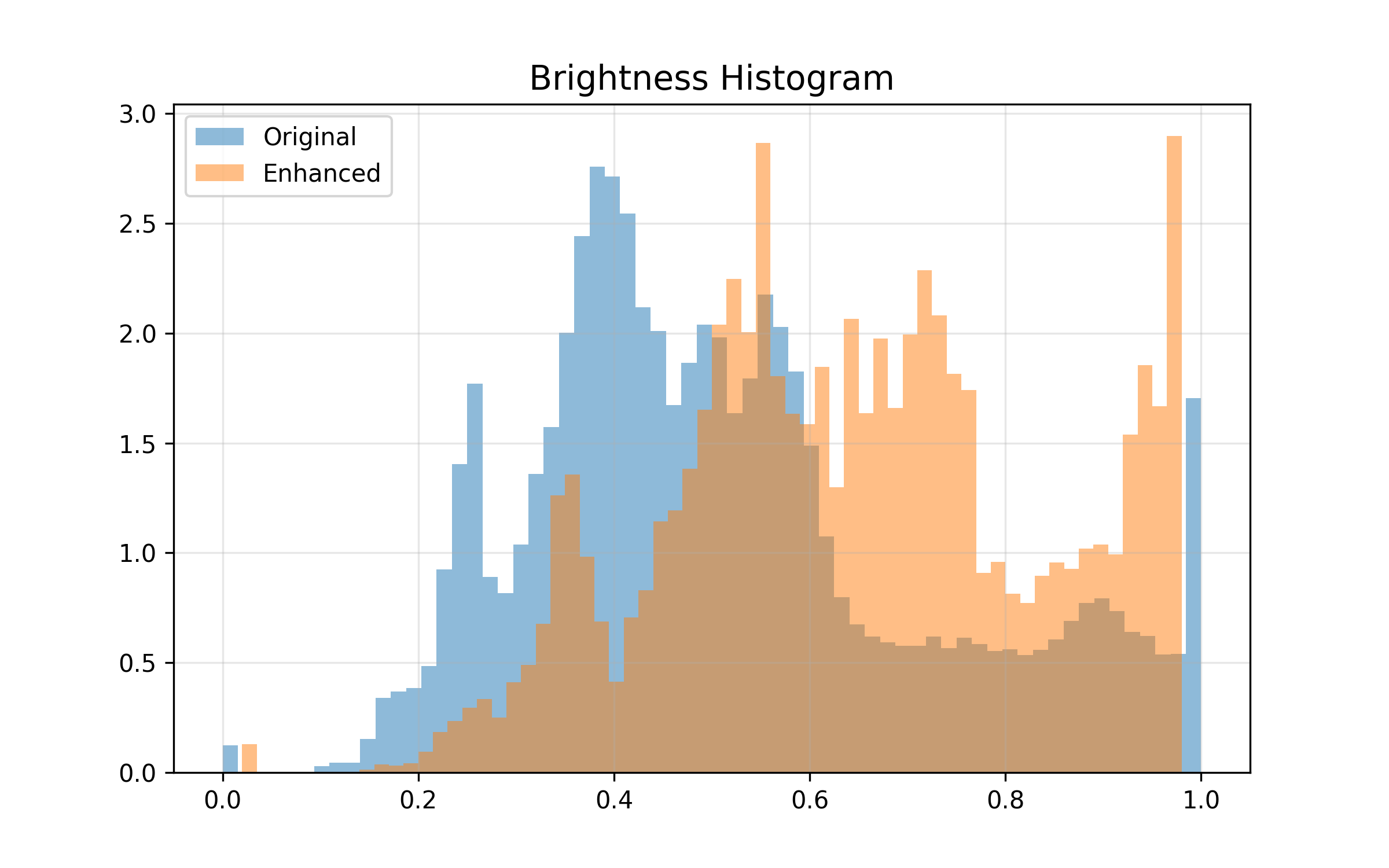} 
			\caption{Brightness histogram}
			\label{fig:brightness_hist}
		\end{subfigure}
		
		\caption{Image enhancement results via the PINN-AFE framework.}
		\label{fig:image_enhance}
	\end{figure}
	
	As shown in Fig. \ref{fig:parrot_compare}, the proposed framework effectively enhances contrast while preserving natural colors and fine details, avoiding over-exposure artifacts common in traditional histogram equalization. The learned mapping curve (Fig. \ref{fig:intensity_curve}) satisfies the strict monotonicity constraint of optimal transport, ensuring no brightness reversal or unnatural artifacts. The brightness histogram (Fig. \ref{fig:brightness_hist}) shows that the method balances the distribution of pixel intensities without introducing the false contouring and blocking artifacts induced by conventional hard histogram mapping. These results validate that the PINN-AFE framework combines mathematical rigor with practical flexibility for image enhancement tasks.
	
	\subsection{Medical Imaging}
	The T1-weighted MRI and FDG-PET images used in this study are obtained from the Harvard Medical School Neuroimaging Primer \citep{johnson_aanlib}, corresponding to a clinical case of mild Alzheimer's disease. The patient presents with typical imaging features: global sulcal widening on T1 MRI, and bilateral hypometabolism in the anterior temporal and posterior parietal cortices on FDG-PET.
	
	In the two-dimensional case for medical image registration, the Monge-Ampère equation takes the following form:
	\begin{equation}
		\det\left(\nabla^2 \psi(\mathbf{x})\right) = \frac{\mu(\mathbf{x})}{\nu\left(\nabla \psi(\mathbf{x})\right)},
	\end{equation}
	where $\nabla^2 \psi(\mathbf{x})$ is the Hessian matrix of the convex potential function $\psi(\mathbf{x})$. The optimal transport map $T(\mathbf{x}) = \nabla \psi(\mathbf{x})$ is inherently a diffeomorphism (invertible, topology-preserving transformation) with strictly positive Jacobian determinant, which is the core property required for artifact-free medical image registration. This framework is equally applicable to T2-weighted MRI registration, as it relies on matching intensity probability distributions rather than specific tissue contrast mechanisms.
	
	For the multimodal image registration task:
	\begin{itemize}[itemsep=0pt, parsep=0pt, topsep=0pt]
		\item $\mu(\mathbf{x})$ is the intensity probability density function of the fixed image
		\item $\nu(\mathbf{y})$ is the intensity probability density function of the moving image
		\item $T(\mathbf{x})$ is the deformation field mapping function to be learned
	\end{itemize}
	
\begin{figure}[H]
		\centering
		\includegraphics[width=0.6\textwidth]{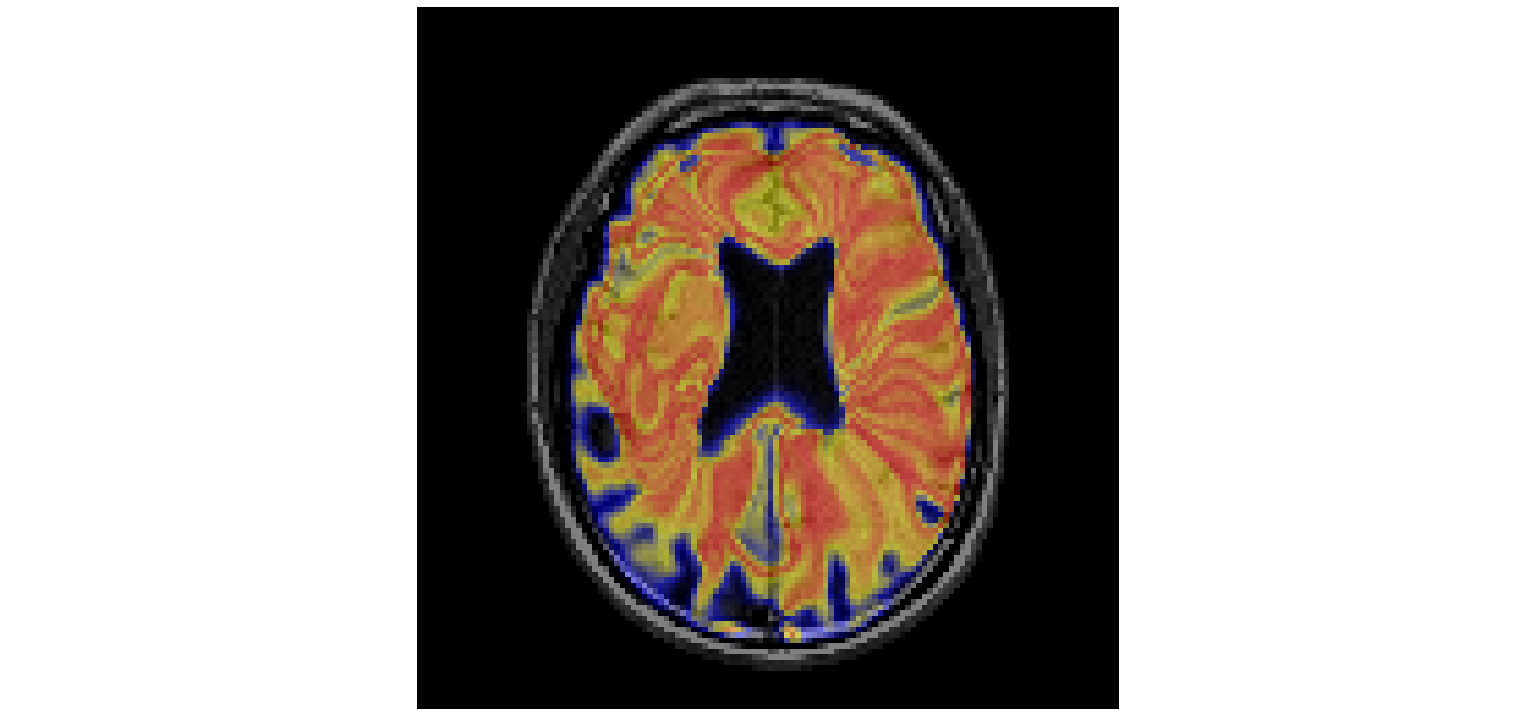}
		\caption{Final overlay}
		\label{fig:pet_overlay}
	\end{figure}

	Table \ref{tab:pet_quantitative} summarizes the quantitative metrics for the clinical dataset after 2000 training epochs, and Figure \ref{fig:pet_overlay} shows the final overlay. The results confirm high registration accuracy including DSC of 0.8635 and Jaccard of 0.7597, good boundary alignment with HD95 of 7.81 px superior to SPM12's 7.89 px, and physically plausible deformation with Fold Ratio of 0 and Jacobian Mean of 0.9969, ensuring the hypometabolic regions in FDG PET are accurately aligned with the anatomical structures in T1 MRI.
	
	\begin{table}[H]
		\centering
		\caption{Quantitative T1-FDG PET Registration Metrics}
		\label{tab:pet_quantitative}
		\resizebox{0.35\textwidth}{!}{
			\begin{tabular}{lc}
				\toprule
				\textbf{Metric} & \textbf{Value} \\
				\midrule
				DSC  & 0.8635 \\
				Jaccard  & 0.7597 \\
				HD95  & 7.81 \\
				Fold Ratio  & 0.0000 \\
				Jacobian Mean  & 0.9969 \\
				Flow Magnitude Mean  & 0.1110 \\
				Smoothness  & 0.0002 \\
				\bottomrule
			\end{tabular}
		}
	\end{table}
	
	\begin{figure}[H]
		\centering
		\setlength{\tabcolsep}{4pt} 
		\begin{tabular}{cc}
			\includegraphics[width=0.45\textwidth]{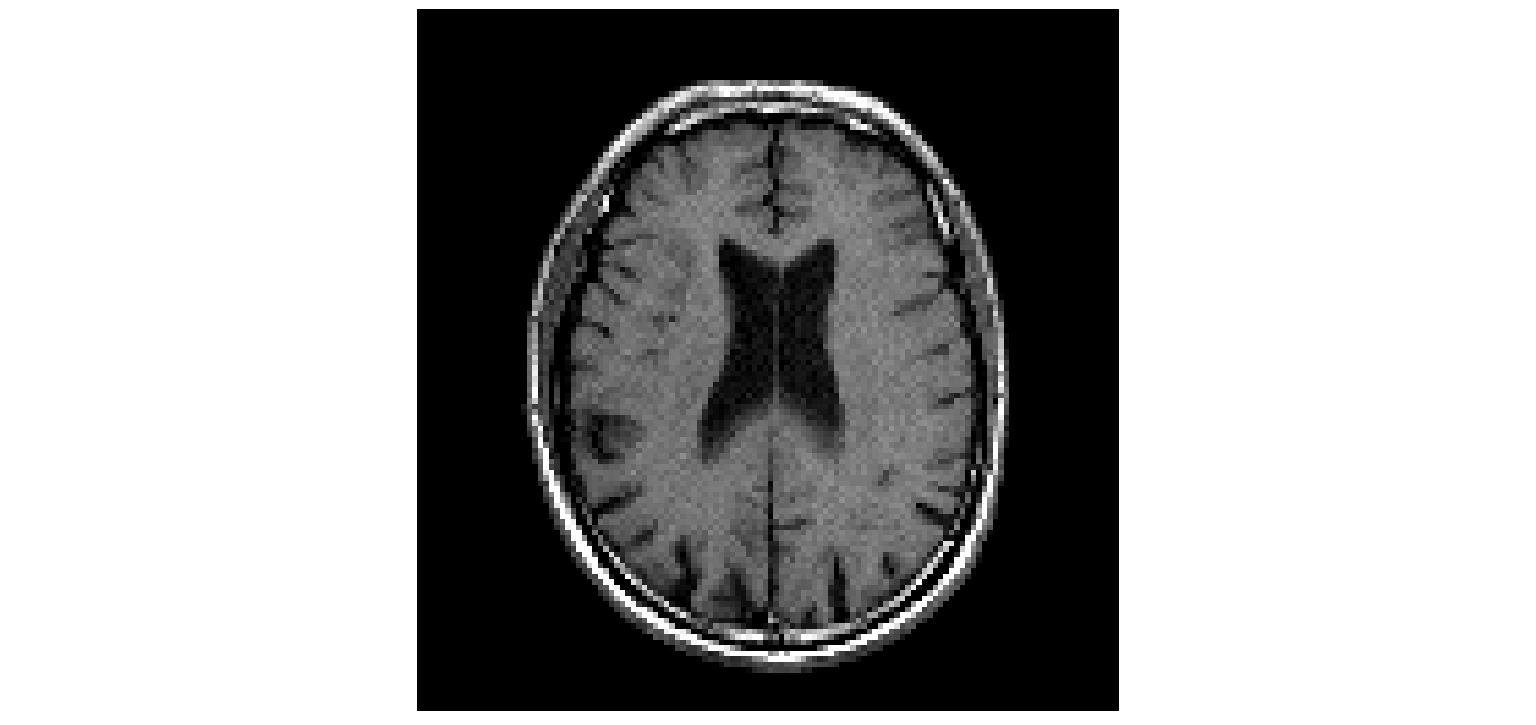} &
			\includegraphics[width=0.45\textwidth]{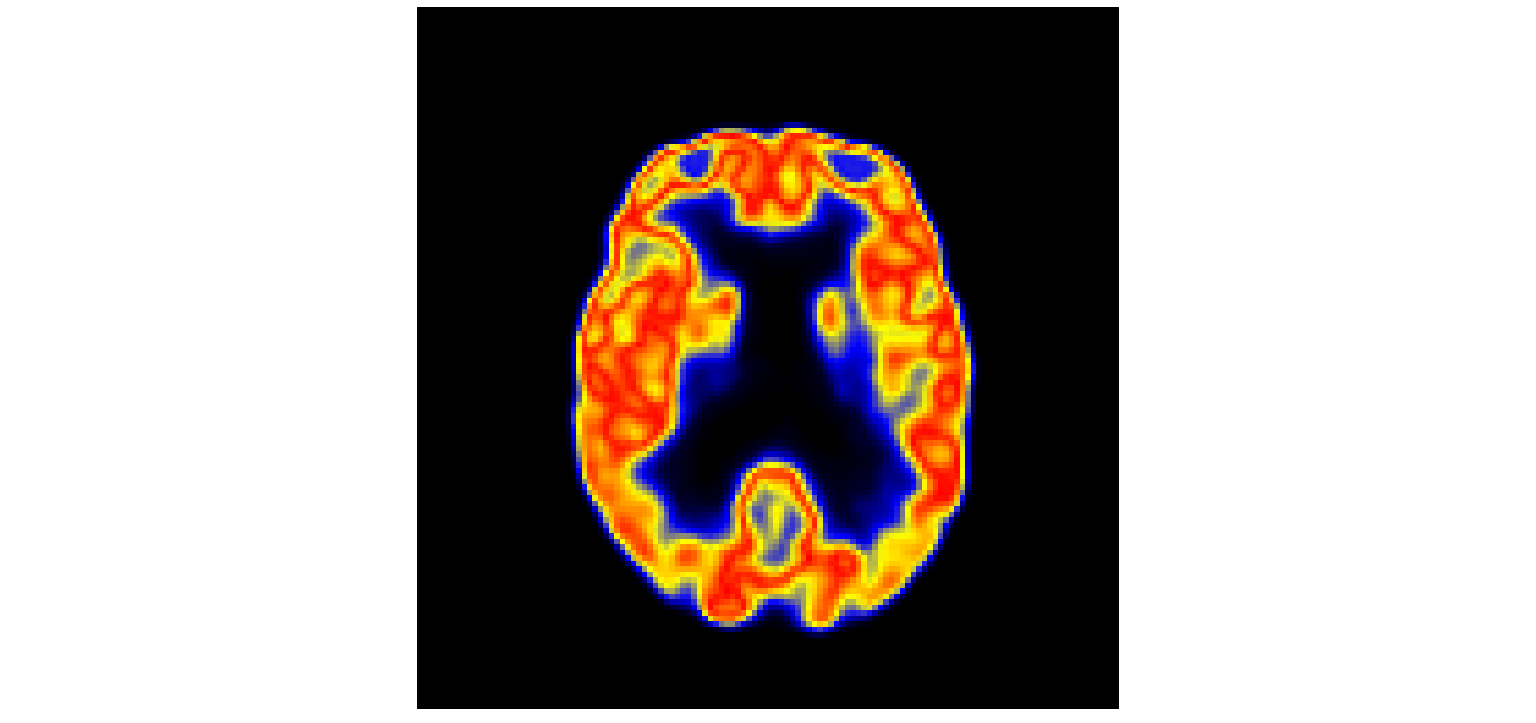} \\
			\small a. T1 Image & \small b. Original FDG Image \\[6pt]
			
			\includegraphics[width=0.45\textwidth]{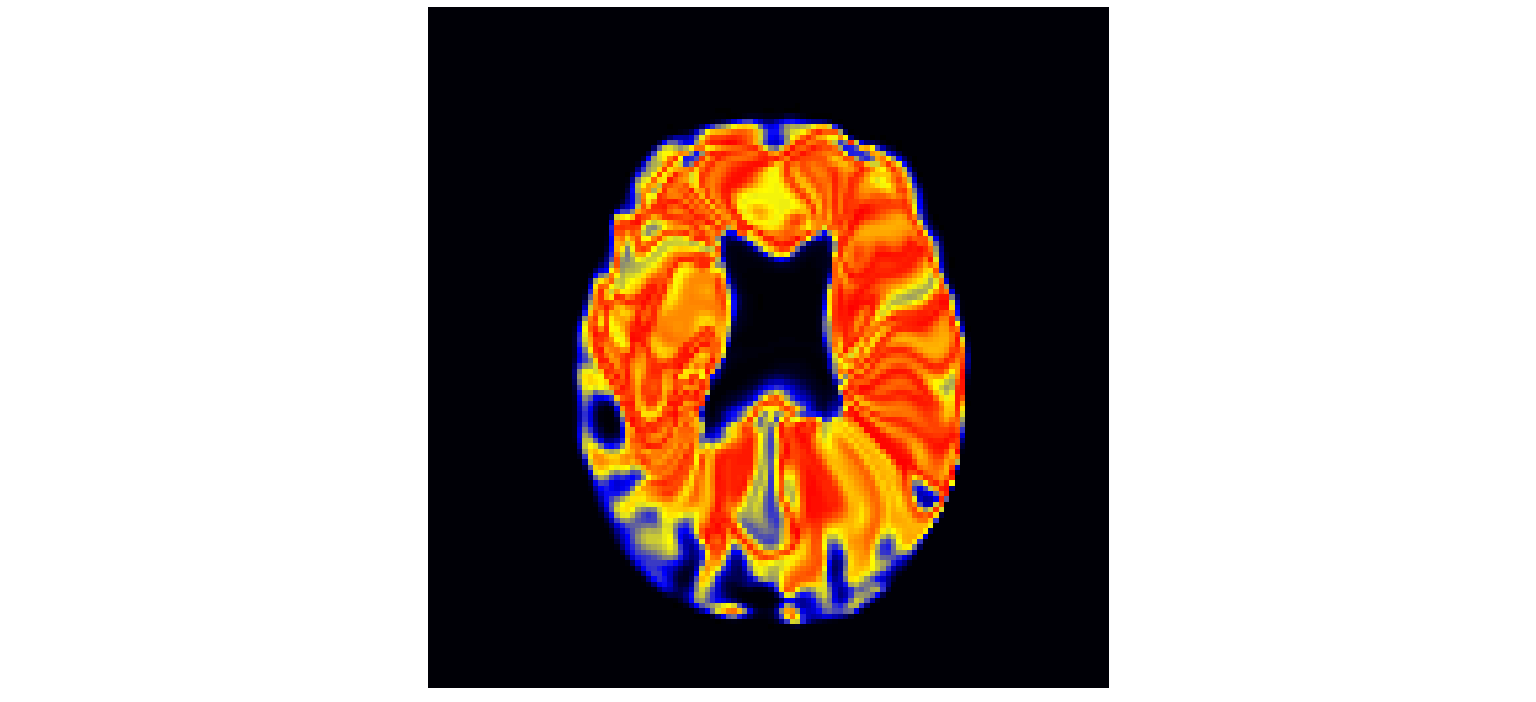} &
			\includegraphics[width=0.45\textwidth]{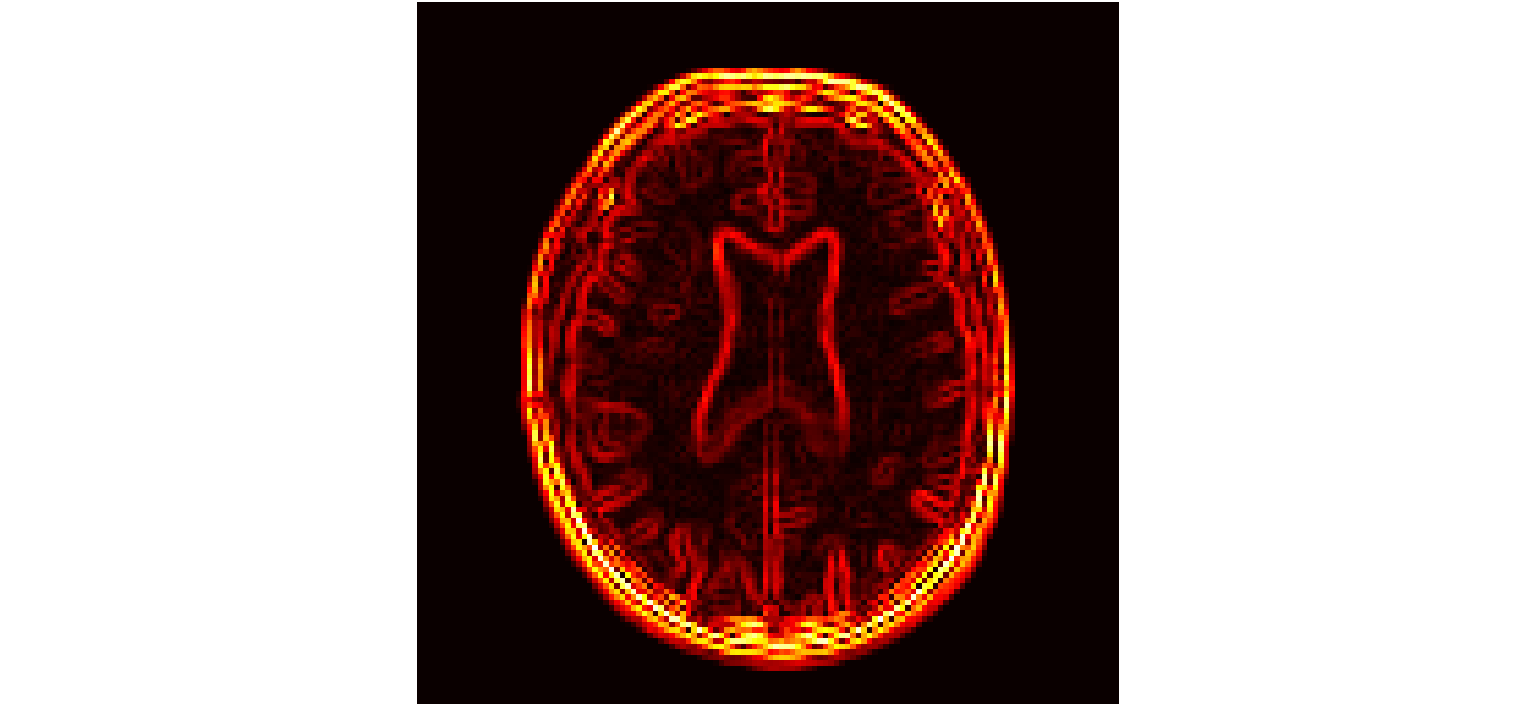} \\
			\small c. Warped FDG Image & \small d. T1 Boundary Map \\[6pt]
			
			\includegraphics[width=0.45\textwidth]{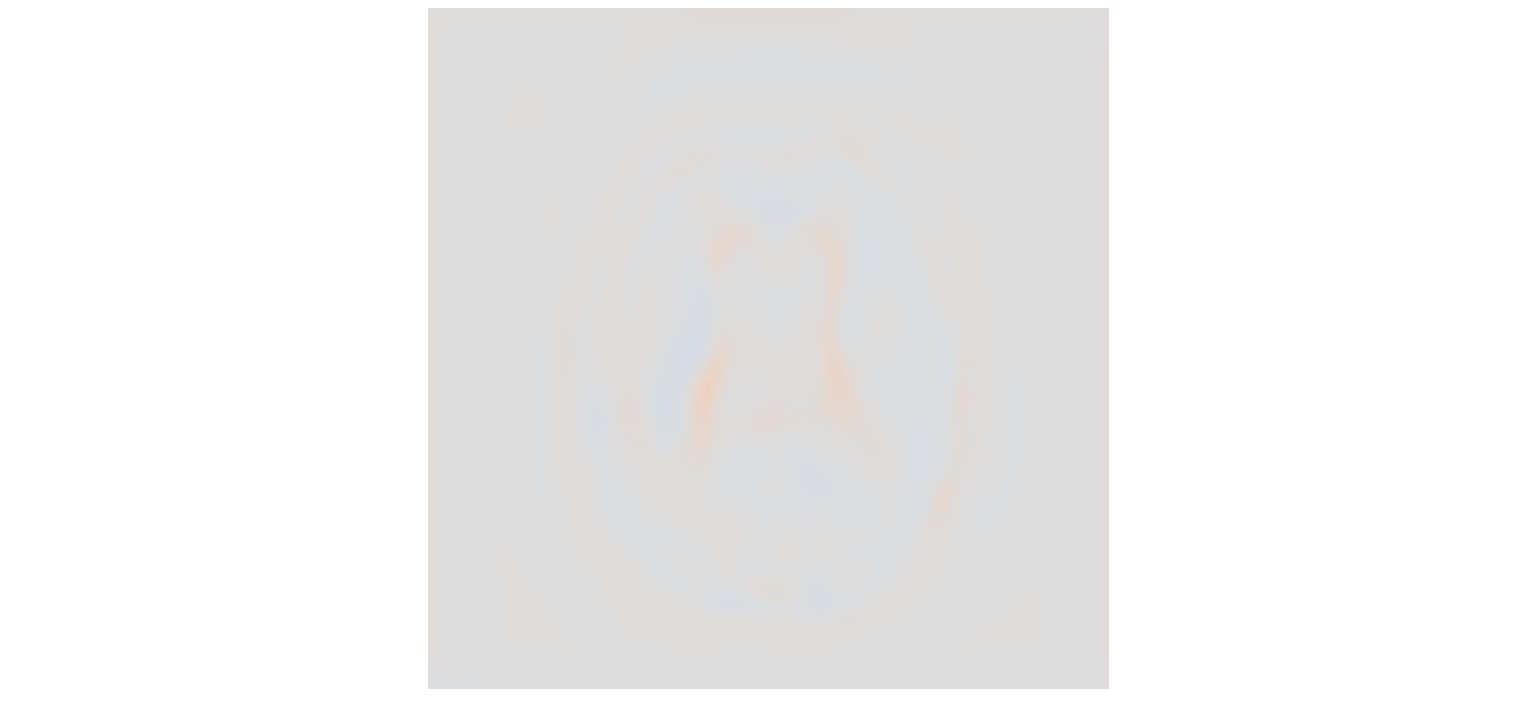} &
			\includegraphics[width=0.45\textwidth]{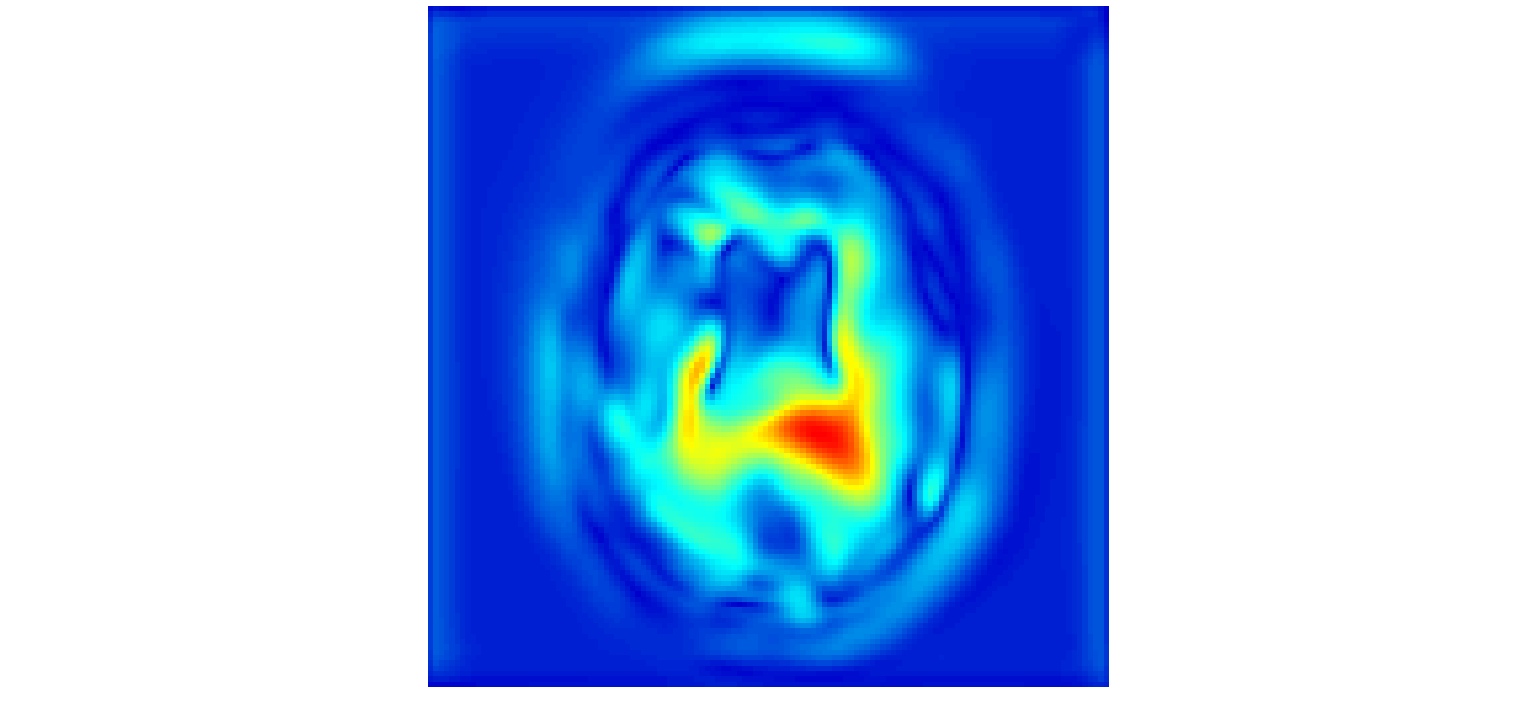} \\
			\small e. Jacobian Determinant & \small f. Flow Magnitude \\
		\end{tabular}
		\caption{Qualitative visualization of T1-FDG PET registration}
		\label{fig:pet_vis}
	\end{figure}
	
	\section{Conclusion}\label{sec:conclusion}
	This paper presents PINN-AFE, a physics-informed neural network with attention feature expansion for the fully nonlinear Monge–Ampère equation.
	
	Theoretical analysis establishes the framework's soundness: error is decomposed into controllable approximation, feature expansion, and optimization components; attention feature expansion reduces residual bounds by a factor of $K$ (number of heads) and lowers parameter/sample complexity from $\mathcal{O}(\varepsilon^{-d/\alpha})$ to $\mathcal{O}(\varepsilon^{-m/\alpha})$ ($m\ll d$); and dynamic weights guarantee a strictly larger effective convergence rate than uniform weighting. Numerical experiments on smooth 2D/3D and singular benchmarks show PINN-AFE achieves MAEs of $10^{-6}$, $1.2\times10^{-5}$, and below $9\times10^{-4}$ respectively, outperforming standard and ICNN-based PINNs by 1–2 orders of magnitude.
	
	Beyond PDE solving, the ICNN-enforced monotonicity prior enables color image enhancement with balanced contrast, brightness, and shadow preservation. For clinical medical image registration, PINN-AFE formulates T1/T2 MRI alignment and T1-FDG PET fusion as Monge–Ampère optimal transport problems, achieving DSC up to 0.8635 and HD95 of 7.81\,px with fold-free deformations, outperforming ANTs and SPM. These results bridge rigorous PDE theory and clinical demands, providing a general methodology for convexity-constrained, attention-driven PINN solvers and deformable registration tasks.
	
	\bibliographystyle{plainnat}
	
	\bibliography{ref}
	
\end{document}